\documentclass[11pt]{article}
\usepackage[LGR,T1]{fontenc}
\usepackage[greek,english]{babel}
\RequirePackage{amsthm,amsmath,natbib}
\RequirePackage[colorlinks,citecolor=blue,urlcolor=blue]{hyperref}
\RequirePackage{hypernat}
\usepackage{amssymb,amsmath,latexsym,color,textcomp,marvosym}
\usepackage[english]{babel}
\usepackage[T1]{fontenc}
\usepackage{tabulary,ulem}
\usepackage{enumerate}
\usepackage{graphicx}
\usepackage{url}

\newcommand{\bepsilon}{\boldsymbol{\epsilon}}

\newcommand{\balpha}{\boldsymbol{\alpha}}
\newcommand{\btheta}{\boldsymbol{\theta}}

\newcommand{\bbeta}{\boldsymbol{\beta}}

\newcommand{\bb}{\boldsymbol{b}}
\newcommand{\bsigma}{\boldsymbol{\sigma}}

\newcommand{\be}{\boldsymbol{e}}
\newcommand{\bu}{\boldsymbol{u}}
\newcommand{\bz}{\boldsymbol{z}}
\newcommand{\bJ}{\boldsymbol{J}}
\newcommand{\bI}{\boldsymbol{I}}

\newcommand{\bP}{\boldsymbol{P}}
\newcommand{\bA}{\boldsymbol{A}}

\newcommand{\bM}{\boldsymbol{M}}
\newcommand{\bN}{\boldsymbol{N}}
\newcommand{\bQ}{\boldsymbol{Q}}

\newcommand{\bC}{\boldsymbol{C}}

\newcommand{\bzero}{\mathbf{0}}
\def\zak{\null\hfill{$\Box$}\par\vspace*{0.2cm}}

\def \argmin{\mathop{\hbox{\rm arg min}}}

\def\1g{1\hskip -3pt \mbox{l}}

\def \Sum {\displaystyle \sum }

\newtheorem{theo}{Theorem}
\newtheorem{lem}{Lemma}
\newtheorem{prop}{Proposition}

\newtheorem{rem}{Remark}
\newtheorem{exemple}{Example}
\numberwithin{equation}{section}
\numberwithin{theo}{section}
\numberwithin{lem}{section}
\numberwithin{cor}{section}
\numberwithin{exemple}{section}
\numberwithin{rem}{section}
\numberwithin{prop}{section}

\renewcommand{\baselinestretch}{1.5}
\begin{document}
\selectlanguage{english}

\title{{
GARCH models without positivity constraints: Exponential or Log
GARCH? }}
\author{
{\sc Christian Francq\footnote{CREST and University Lille 3
(EQUIPPE), BP 60149, 59653 Villeneuve d'Ascq cedex, France. E-Mail:
christian.francq@univ-lille3.fr}, Olivier Wintenberger
\footnote{CEREMADE (University Paris-dauphine) and CREST, Place du
Maréchal de Lattre de Tassigny, 75775 PARIS Cedex 16, France.
 E-mail: wintenberger@ceremade.dauphine.fr} and Jean-Michel
Zakoïan\footnote{Corresponding author: Jean-Michel Zakoïan, EQUIPPE
(University Lille 3) and CREST, 15 boulevard Gabriel Péri, 92245
Malakoff Cedex, France. E-mail: zakoian@ensae.fr, Phone number:
33.1.41.17.77.25.} }}

\selectlanguage{english}

\date{} \maketitle
\begin{quote}
\begin{center}
\textbf{Abstract}
\end{center}
{\small \hspace{1em}  This paper provides a probabilistic and
statistical comparison of the log-GARCH and EGARCH models, which
both rely on multiplicative volatility dynamics without positivity
constraints.
We compare the main probabilistic properties (strict stationarity,
existence of moments, tails) of the EGARCH model, which are already
known, with those of an asymmetric version
of the log-GARCH. 
The quasi-maximum likelihood estimation of the log-GARCH parameters
is shown to be strongly consistent and asymptotically normal.
Similar estimation results are only available for the EGARCH(1,1)
model, and under much stronger assumptions. The comparison is
pursued via simulation experiments and estimation on real data. }
\end{quote}


\medskip


\newpage
\section{Preliminaries}

Since their introduction by Engle (1982) and Bollerslev (1986),
GARCH models have attracted much attention and have been widely
investigated in the literature. 
 Many extensions have been suggested
and, among them, the EGARCH (Exponential GARCH) introduced and
studied by Nelson (1991) is very popular. In this model, the
log-volatility is expressed as a linear combination of its past
values and past values of the positive and negative parts of the
innovations. Two main reasons for the success of this formulation
are that (i) it allows for asymmetries in volatility (the so-called
leverage effect: negative shocks tend to have more impact on
volatility than positive shocks of the same magnitude), and (ii) it
does not impose any positivity restrictions on the volatility coefficients.

Another class of GARCH-type models, which received less attention,
seems to share the same characteristics. The log-GARCH(p,q) model
has been introduced, in slightly different forms, by Geweke (1986),
Pantula (1986) and Milhøj (1987). For more recent works on this
class of models, the reader is referred to Sucarrat and Escribano
(2010) and the references therein. The (asymmetric) log-GARCH$(p,q)$
model takes the form
\begin{equation}\label{logGARCH}
\left\{\begin{array}{lll}
\epsilon_t&=&\sigma_t\eta_t,\quad \\
\log \sigma_t^2&=&\omega+\sum_{i=1}^q\left(\alpha_{i+}1_{\{\epsilon_{t-i}>0\}}+\alpha_{i-}1_{\{\epsilon_{t-i}<0\}}\right)\log\epsilon_{t-i}^2\\
&&+\sum_{j=1}^p\beta_j\log \sigma_{t-j}^2
\end{array}\right.
\end{equation} where $\sigma_t>0$ and $(\eta_t)$ is a
sequence of independent and identically distributed (iid) variables
such that $E\eta_0=0$ and $E\eta_0^2=1$. The usual symmetric
log-GARCH corresponds to the case $\balpha_{+}=\balpha_{-}$, with
$\balpha_{+}=(\alpha_{1+},\dots,\alpha_{q+})$ and
$\balpha_{-}=(\alpha_{1-},\dots,\alpha_{q-})$.

Interesting features of the log-GARCH specification
are the following.

{\bf (a) Absence of positivity constraints.}
 An advantage of
modeling the log-volatility rather than the volatility is that the
vector $\btheta=(\omega,\balpha_+,\balpha_-,\bbeta)$ with
$\bbeta=(\beta_{1},\dots,\beta_{p})$ is not a priori subject to
positivity constraints\footnote{However, some desirable properties
may determine the sign of  coefficients. For instance, the present
volatility is generally thought of as an increasing function of its
past values, which entails $\beta_j>0.$ The difference with standard
GARCH models is that such constraints are not required for the
existence of the process and, thus, do not complicate estimation
procedures.}. This property seems particularly appealing when
exogenous variables are included in the volatility specification
(see Sucarrat and Escribano, 2012).

{\bf (b) Asymmetries.} Except when $\alpha_{i+}=\alpha_{i-}$ for all
$i$, positive and negative past values of $\epsilon_t$ have
different impact on the current log-volatility, hence on the current
volatility. However, given that $\log\epsilon_{t-i}^2$ can be
positive or negative, the usual leverage effect does not have a
simple characterization, like $\alpha_{i+}<\alpha_{i-}$ say. Other
asymmetries could be introduced, for instance by replacing $\omega$
by
$\sum_{i=1}^q\omega_{i+}1_{\{\epsilon_{t-i}>0\}}+\omega_{i-}1_{\{\epsilon_{t-i}<0\}}$.
The model would thus be stable by scaling, which is not the case of
Model (\ref{logGARCH}) except in the symmetric case.

{\bf (c) The volatility is not bounded below.} Contrary to standard
GARCH models and most of their extensions, there is no minimum value
for the volatility. The existence of such a bound can be problematic
because, for instance in a GARCH(1,1), the minimum value is
determined by the intercept $\omega$. On the other hand, the
unconditional variance is proportional to $\omega$. Log-volatility
models allow to disentangle these two properties (minimum value and
expected value of the volatility).

{\bf (d) Small values can have persistent effects on volatility.} In
usual GARCH models, a large value (in modulus) of the volatility
will be followed by other large values (through the coefficient
$\beta$ in the GARCH(1,1), with standard notation). A sudden rise of
returns (in module) will also be followed by large volatility values
if the coefficient $\alpha$ is not too small. We thus have
persistence of large returns and volatility. But small returns (in
module) and small volatilities are not persistent. In a period of
large volatility, a sudden drop of the return due to a small
innovation, will not much alter the subsequent volatilities (because
$\beta$ is close to 1 in general). By contrast, as will be
illustrated in the sequel, the log-GARCH provides persistence of
large {\it and} small values.

{\bf (e) Power-invariance of the volatility specification.}  An
interesting potential property of time series models is their
stability with respect to certain transformations of the
observations. Contemporaneous aggregation and temporal aggregation
of GARCH models have, in particular, been studied by several authors
(see  Drost and Nijman (1993)). On the other hand, the choice of a
power-transformation is an issue for the volatility specification.
For instance, the volatility can be expressed in terms of past
squared values (as in the usual GARCH) or in terms of past absolute
values (as in the symmetric TGARCH) but such specifications are
incompatible. On the contrary, any power transformation
$|\sigma_t|^s$ (for $s\ne 0$) of a log-GARCH volatility has a
log-GARCH form (with the same coefficients in $\btheta$, except the
intercept $\omega$ which is multiplied by $s/2$).

\begin{figure}[!h]
\includegraphics[width=8.5cm]{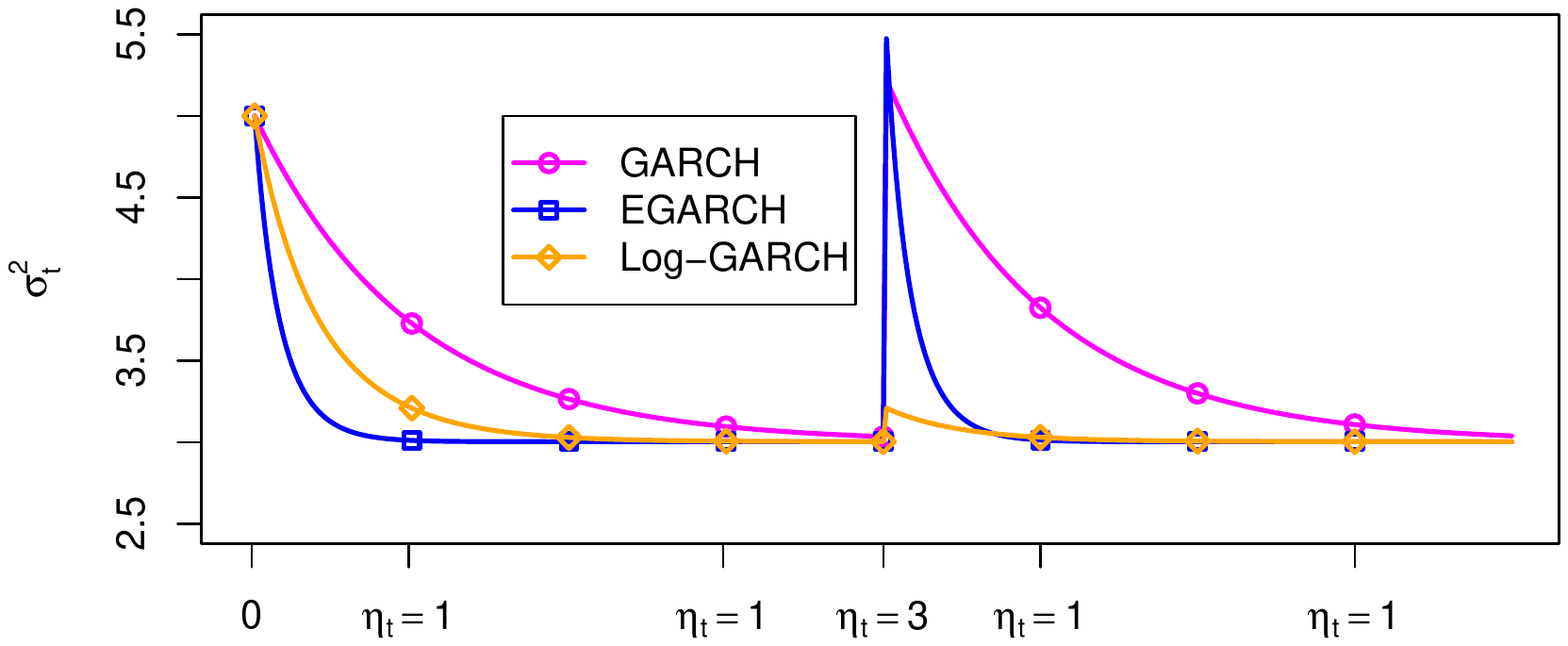}
\includegraphics[width=8.5cm]{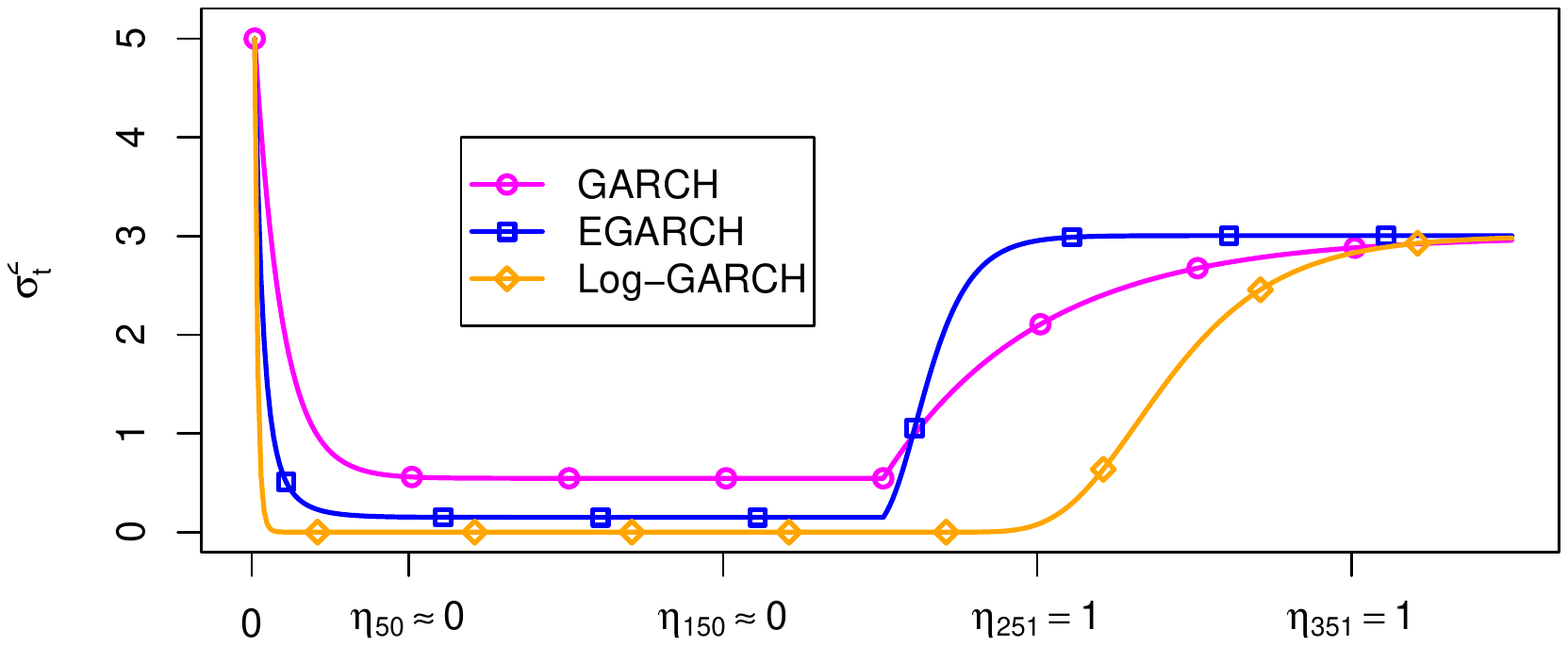}
\begin{center}{
\includegraphics[width=8.5cm]{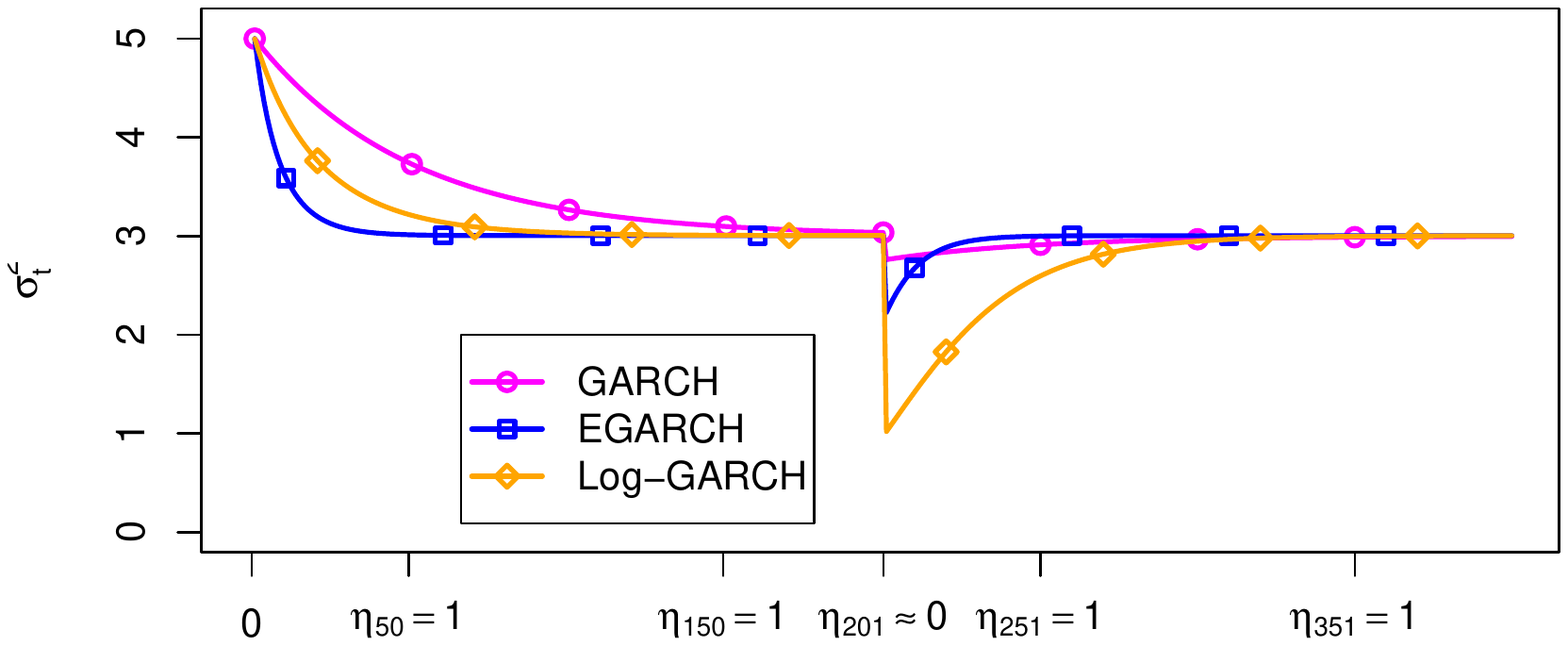}}\end{center}\vspace{-1cm}
\caption{Curves of the impact of shocks on volatility. 
}
\label{impcur}
\end{figure}

The log-GARCH model has apparent similarities with the
EGARCH($p,\ell$) model defined by
\begin{equation}\label{EGARCH}
\left\{\begin{array}{lll}
\epsilon_t&=&\sigma_t\eta_t,\quad \\
\log \sigma_t^2&=&\omega+
\sum_{j=1}^p\beta_{j}\log \sigma_{t-j}^2
+ \sum_{k=1}^{\ell} \gamma_{k} \eta_{t-k}+ \delta_{k} |\eta_{t-k}|,
\end{array}\right.
\end{equation}
under the same assumptions on the sequence $(\eta_t)$ as in Model
(\ref{logGARCH}). These models have in common the above properties
\textbf{(a)}, \textbf{(b)}, \textbf{(c)} and \textbf{(e)}.
Concerning the property in \textbf{(d)}, and more generally the
impact of shocks on the volatility dynamics,  Figure \ref{impcur}
illustrates the differences between the two models (and also with
the standard GARCH). The  coefficients of the GARCH(1,1) and the
symmetric EGARCH(1,1) and log-GARCH(1,1) models have been chosen to
ensure the same long-term variances when the squared innovations are
equal to 1. Starting from the same initial value $\sigma_0^2$, we
analyze the effect of successive shocks $\eta_t, t\geq 1.$ The
top-left graph shows that a sudden large shock, in the middle of the
sample, has a (relatively) small impact on the log-GARCH, a large
but transitory effect on the EGARCH,
 and a large and very persistent effect on the classical GARCH volatility.
 The top-right graph shows the effect of a sequence of tiny innovations, $\eta_t \approx 0$ for $t\leq 200$: for the log-GARCH, contrary to the GARCH and EGARCH,
 the effect is persistent. The bottom graph shows that even one tiny innovation causes this persistence of small volatilities for the
 log-GARCH,
contrary to the EGARCH and
GARCH volatilities. 

This article provides a probability and statistical study of the
log-GARCH, together with a comparison with the  EGARCH. While the
stationarity properties of the EGARCH are well-known, those of the
asymmetric log-GARCH$(p,q)$ model (\ref{logGARCH}) have not yet been
established, to our knowledge. 
As for the quasi-maximum likelihood estimator (QMLE), the
consistency and asymptotic normality have only been proved in
particular cases and under cumbersome assumptions for the EGARCH,
but, except in the log-ARCH case by Kristensen and Rahbek (2009),
have not yet been established for the log-GARCH.
Finally, it seems important to compare the two classes of models on
typical financial series. The distinctive features of the two models
may render one or the other formulation more adequate for certain types of
series.

 The remainder of the paper is organized as follows.
Section \ref{sec2} studies the existence of  a solution to Model
(\ref{logGARCH}). Conditions for the existence of log-moments are
derived,  and we characterize the leverage effect. Section
\ref{secdep} is devoted to the tail properties of the solution. In
Section \ref{sec3}, the strong consistency and the asymptotic
normality of the QMLE are established under mild conditions. Section
\ref{sec7} presents some numerical applications on simulated and
real data. Proofs are collected in Section \ref{secProofs}. Section
\ref{sec8} concludes.

\section{Stationarity, moments and asymmetries of the log-GARCH}
\label{sec2}

We start by studying the existence of solutions to Model
(\ref{logGARCH}).

\subsection{Strict stationarity}
Let $\bzero_k$ denote a $k$-dimensional vector of zeroes, and let
$\bI_k$ denote the $k$-dimensional identity matrix. Introducing the
vectors
\begin{eqnarray*}
\bepsilon_{t,q}^+&=&(1_{\{\epsilon_{t}>0\}}\log \epsilon_t^2,\dots,1_{\{\epsilon_{t-q+1}>0\}}\log \epsilon_{t-q+1}^2)'\in \mathbb{R}^q,\\
\bepsilon_{t,q}^-&=&(1_{\{\epsilon_{t}<0\}}\log \epsilon_t^2,\dots,1_{\{\epsilon_{t-q+1}<0\}}\log \epsilon_{t-q+1}^2)'\in \mathbb{R}^q,\\
\bz_t&=&(\bepsilon_{t,q}^+,\bepsilon_{t,q}^-,\log \sigma_t^2,\dots,\log \sigma_{t-p+1}^2)'\in \mathbb{R}^{2q+p},\\
\bb_t&=&\left((\omega+\log\eta_t^2)1_{\{\eta_{t}>0\}},\bzero'_{q-1},(\omega+\log\eta_t^2)1_{\{\eta_{t}<0\}},\bzero'_{q-1},\omega,\bzero'_{p-1}\right)'\in \mathbb{R}^{2q+p},
\end{eqnarray*} and the matrix
\begin{equation}
\label{A_tmulti}
\bC_t=\left(\begin{array}{ccc}
1_{\{\eta_{t}>0\}}\balpha_+&1_{\{\eta_{t}>0\}}\balpha_-&1_{\{\eta_{t}>0\}}\bbeta  \\
\begin{array}{cc}\bI_{q-1} &\bzero_{q-1}\end{array}  &\bzero_{(q-1)\times q}&\bzero_{(q-1)\times p}\\
1_{\{\eta_{t}<0\}}\balpha_+&1_{\{\eta_{t}<0\}}\balpha_-&1_{\{\eta_{t}<0\}}\bbeta  \\
\bzero_{(q-1)\times q}&\begin{array}{cc}\bI_{q-1} &\bzero_{q-1}\end{array}   &\bzero_{(q-1)\times p}\\
\balpha_+&\balpha_-&\bbeta  \\
\bzero_{(p-1)\times q}   &\bzero_{(p-1)\times q}&\begin{array}{cc}\bI_{p-1} &\bzero_{p-1}\end{array}\\
\end{array}\right),
\end{equation}
we rewrite Model (\ref{logGARCH}) in matrix form as
\begin{equation}\label{vectorform}
\bz_t=\bC_t\bz_{t-1}+\bb_t.\end{equation}
We have implicitly assumed $p>1$ and $q>1$ to write $\bC_t$ and  $\bb_t$, but obvious changes of notation can be employed when $p\leq 1$ or $q\leq 1$.
Let $\gamma ({\bf C})$
be the top Lyapunov exponent of the
sequence ${\bf C}=\{C_t, t\in \mathbb{Z}\}$,
\begin{equation*}
\gamma ({\bf C})=\lim_{t\to \infty} \frac{1}{t} E \left(\log
\|\bC_t\bC_{t-1}\ldots \bC_{1}\|\right) = \inf_{t\geq 1}\;
\frac{1}{t}\; E (\log \|\bC_t\bC_{t-1}\ldots \bC_{1}\|).
\end{equation*}
The choice of the norm is obviously unimportant for the value of the
top Lyapunov exponent. However, in the sequel, the matrix norm will
be assumed to be multiplicative. Bougerol and Picard (1992a) showed
that if an equation of the form (\ref{vectorform}) with iid
coefficients $(\bC_t,\bb_t)$ is  irreducible\footnote{See their
Definition 2.3.}  and if $E\log^+\|\bC_0\|$ and $E\log^+\|\bb_0\|$
are finite, $\gamma ({\bf C})<0$ is the necessary and sufficient
condition for the existence of a stationary solution to
(\ref{vectorform}). Bougerol and Picard (1992b) showed that, for the
univariate GARCH($p,q$) model, there exists a representation of the
form  (\ref{vectorform}) with positive coefficients, and for which
the necessary and sufficient condition for the existence of a
stationary GARCH model is  $\gamma ({\bf C})<0$. The result can be
extended to more general classes of GARCH models (see {\it  e.g.}
Francq and Zakoïan, 2010a). The problem is more delicate with the
log-GARCH because the coefficients of (\ref{vectorform}) are not
constrained to be positive. The following result and Remark
\ref{rk1} below show that $\gamma ({\bf C})<0$ is only sufficient.
The condition is however necessary under the mild additional
assumption that (\ref{vectorform}) is irreducible.
\begin{theo}\label{theostatio}
Assume that $E\log^+|\log\eta_0^2|<\infty$.
A sufficient condition for the
existence of a strictly stationary solution
 to the log-GARCH model (\ref{logGARCH}) is $\gamma ({\bf C})<0$.
When $\gamma ({\bf C})<0$ there exists only one  stationary solution, which is non anticipative  and  ergodic.
\end{theo}
\begin{exemple}[{\bf The log-GARCH(1,1) case}]\label{ex11}
{\rm In the case $p=q=1$, omitting subscripts, we have
$$\bC_t\bC_{t-1}\ldots \bC_{1}=\left(\begin{array}{c}1_{\{\eta_{t}>0\}}\\1_{\{\eta_{t}<0\}}\\ 1\end{array}\right)\left(\begin{array}{ccc}\alpha_+&\alpha_-&\beta\end{array}\right)\prod_{i=1}^{t-1}\left(\alpha_+1_{\{\eta_{i}>0\}}
+\alpha_-1_{\{\eta_{i}<0\}}+\beta\right).$$ Assume that $E\log^+|\log\eta_0^2|<\infty$, which entails $P(\eta_0=0)=0$. Thus,
$$\gamma ({\bf C})=E\log\left|\alpha_+1_{\{\eta_{0}>0\}}
+\alpha_-1_{\{\eta_{0}<0\}}+\beta\right|=\log|\beta+\alpha_+|^a|\beta+\alpha_-|^{1-a},$$
where $a=P(\eta_0>0)$. The condition
$|\alpha_++\beta|^a|\alpha_-+\beta|^{1-a}<1$ thus guarantees the
existence of a stationary solution to the log-GARCH(1,1) model.}
\end{exemple}
\begin{exemple}[{\bf The symmetric case}]
{\rm In the case $\balpha_+=\balpha_-=\balpha$, one can see directly
from (\ref{logGARCH}) that $\log\sigma_t^2$ satisfies an ARMA-type
equation of the form
$$\left\{1-\sum_{i=1}^r\left(\alpha_{i}+\beta_i\right)B^i\right\}\log\sigma_t^2=c+\sum_{i=1}^q\alpha_{i}B^iv_t$$
where $B$ denotes the backshift operator, $v_t=\log\eta_t^2$,
$r=\max{\{p,q\}}$, $\alpha_i=0$ for $i>q$ and $\beta_i=0$ for $i>p$.
This equation is a standard ARMA$(r,q)$ equation under the moment
condition $E(\log\eta_t^2)^2<\infty$, but this assumption is not
needed. It is well known that this equation admits a non degenerated
and non anticipative stationary solution if and only if the roots of
the AR polynomial lie outside the unit circle.

We now show that this condition is equivalent to the condition
$\gamma ({\bf C})<0$ in the case $q=1$. Let $\bP$ be the permutation
matrix obtained by permuting the first and second rows of
$\bI_{2+p}$. Note that
$\bC_t=\bC^+1_{\{\eta_t>0\}}+\bC^-1_{\{\eta_t<0\}}$ with
$\bC^-=\bP\bC^+$. Since $\alpha_+=\alpha_-$, we have
$\bC^+\bP=\bC^+$. Thus $\bC^+\bC^-=\bC^+\bP\bC^+=\bC^+\bC^+$ and
$\|\bC_t\cdots \bC_1\| = \|(\bC^+)^t\|.$ It follows that $\gamma
({\bf C})=\log\rho(\bC^+)$. In view of the companion form of
$\bC^+$, it can be seen that the condition $\rho(\bC^+)<1$ is
equivalent to the condition $z-\sum_{i=1}^r(\alpha_{i}+\beta_i)z^i=0
\Rightarrow |z|>1.$}
\end{exemple}
\begin{rem}[{\bf The condition $\gamma ({\bf C})<0$ is not
necessary}]\label{rk1} {\rm  Assume for instance that $p=q=1$ and
$\alpha_+=\alpha_-=\alpha$. In that case  $\gamma ({\bf C})<0$ is
equivalent to $|\alpha+\beta|<1$. In addition, assume that
$\eta_0^2=1$ a.s. Then, when $\alpha+\beta\neq 1$, there exists a
stationary solution to (\ref{logGARCH}) defined by
$\epsilon_t=\exp(c/2)\eta_t$, with $c=\omega/(1-\alpha-\beta)$.}
\end{rem}

\subsection{Existence of log-moments}
It is well known that for GARCH-type models, the strict stationarity condition entails the existence of a moment of order $s>0$ for $|\epsilon_t|$. The following Lemma shows that this is also the case for $|\log\epsilon_t^2|$ in the log-GARCH model, when the condition $E\log^+|\log\eta_0^2|<\infty$ of Theorem~\ref{theostatio} is slightly reinforced.
\begin{prop}[{\bf Existence of a fractional log-moment}]\label{fracmoment}
Assume that $\gamma ({\bf C})<0$ and that $E|\log\eta_0^2|^{s_0}<\infty$ for some $s_0>0$. Let $\epsilon_t$ be the strict stationary solution of (\ref{logGARCH}). There exists $s>0$ such that $E|\log\epsilon_t^2|^s<\infty$ and $E|\log\sigma_t^2|^s<\infty$.
\end{prop}

In order to give conditions for the existence of higher-order
moments, we  introduce some additional notation. Let $\be_i$ be the
$i$-th column of $\bI_r$, let $\bsigma_{t,r}=(\log \sigma_t^2,\dots,
\log \sigma_{t-r+1}^2)'$, $r=\max{\{p,q\}}$, and let the companion
matrix
\begin{align}
\label{MatA}
\bA_t&=
\left(\begin{array}{cccc}
\mu_1(\eta_{t-1})&\dots& \mu_{r-1}(\eta_{t-r+1})& \mu_r(\eta_{t-r})\\
\multicolumn{3}{c}{\bI_{r-1}} & \bzero_{r-1}
\end{array}\right),&
\end{align}
where
$\mu_i(\eta_{t})=\alpha_{i+}1_{\{\eta_{t}>0\}}+\alpha_{i-}1_{\{\eta_{t}<0\}}+\beta_i$
with the convention $\alpha_{i+}=\alpha_{i-}=0$ for $i>p$ and
$\beta_{i}=0$ for $i>q$. We have the Markovian representation
\begin{equation}
\label{Markovian}
\bsigma_{t,r}=\bA_t\bsigma_{t-1,r}+\bu_t,
\end{equation}
where
$\bu_t=u_t\be_1$,
with
$$u_t=\omega+\sum_{i=1}^q\left(\alpha_{i+}1_{\{\eta_{t-i}>0\}}+\alpha_{i-}1_{\{\eta_{t-i}<0\}}\right)\log\eta_{t-i}^2.$$
The sequence of matrices $(\bA_t)$ is dependent, which makes
(\ref{Markovian}) more difficult to handle than (\ref{vectorform}).
On the other hand, the size of the  matrices $\bA_t$ is smaller than
that of $\bC_t$ ($r$ instead of $2q+p$) and, as we will see, the
log-moment conditions obtained with (\ref{Markovian}) can be sharper
than with (\ref{vectorform}).

Before deriving such log-moment conditions, we need some additional
notation.
 The Kronecker matrix product is denoted by
$\otimes$, and the spectral radius of a square matrix $\bM$ is
denoted by $\rho(\bM)$. Let $M^{\otimes m}=M\otimes \ldots \otimes
M$. For any (random) vector or matrix $\bM$, let
$\mathrm{Abs}({\bM})$ be the matrix, of same size as $\bM$, whose
elements are the absolute values of the corresponding elements of
$\bM$. For any sequence of identically distributed random matrices
matrices $(\bM_t)$ and for any integer $m$,  let
$\bM^{(m)}=E[\{\mathrm{Abs}({\bM}_{1})\}^{\otimes m}]$.
\begin{prop}[{\bf Existence of $m$-order log-moments}]\label{logmoment}
Let $m$ be a positive integer. Assume that $\gamma ({\bf C})<0$ and that $E|\log\eta_0^2|^{m}<\infty$. 
\begin{itemize}
\item If $m=1$ or $r=1$, then $\rho(\bA^{(m)})<1$ implies that
the strict stationary solution of (\ref{logGARCH}) is such that
$E|\log\epsilon_t^2|^m<\infty$ and $E|\log\sigma_t^2|^m<\infty$.
\item If $\rho(\bC^{(m)})<1$, then $E|\log\epsilon_t^2|^m<\infty$ and
$E|\log\sigma_t^2|^m<\infty$.
\end{itemize}
\end{prop}
Note that the conditions $\rho(\bA^{(m)})<1$ and $\rho(\bC^{(m)})<1$
are similar to those obtained by Ling and McAleer (2002) for the
existence of moments of standard GARCH models.
\begin{exemple}[{\bf Log-GARCH(1,1) continued}]
\label{examplelog-garch} In the case $p=q=1$,  we have
$$\bA_t=\alpha_+1_{\{\eta_{t-1}>0\}}+\alpha_-1_{\{\eta_{t-1}<0\}}+\beta \quad\mbox{and}\quad \bA^{(m)}=E\left(|\bA_1|\right)^m.$$
The conditions $E|\log\eta_0^2|^m<\infty$ and, with $a=P(\eta_0>0)$,
$$a\left|\alpha_++\beta\right|^m+(1-a)\left|\alpha_-+\beta\right|^{m}<1$$ thus entail $E|\log\epsilon_t^2|^m<\infty$ for the log-GARCH(1,1) model.
Note that the condition $\rho(\bC^{(m)})<1$ takes the (more binding)
form
$$a\left(\left|\alpha_+\right|+\left|\beta\right|\right)^m+(1-a)\left(\left|\alpha_-\right|+\left|\beta\right|\right)^{m}<1$$
\end{exemple}
Now we study the existence of any log-moment. Let
$\bA^{(\infty)}=\mathrm{ess}\sup \mathrm{Abs}({\bA}_{1})$ be the
essential supremum of $\mathrm{Abs}({\bA}_{1}) $ term by term.
\begin{prop}[{\bf Existence of log-moments at any order}]\label{anylogmoment}
Suppose that $\gamma ({\bf C})<0$ and
\begin{equation}\label{allmombis}
\rho(\bA^{(\infty)})<1
\quad \mbox{or, equivalently, } \quad
\sum_{i=1}^r \left|\alpha_{i+}+\beta_i\right|\vee
\left|\alpha_{i-}+\beta_i\right|<1.
\end{equation}
Then $E|\log\epsilon_t^2|^m<\infty$ at any order $m$ such that
$E|\log\eta_0^2|^{m}<\infty$.
\end{prop}

\subsection{Leverage effect}
A well-known stylized fact of financial markets is that negative
shocks on the returns impact  future volatilities more importantly
than positive shocks of the same magnitude. In the log-GARCH(1,1)
model with $\alpha_->\max\{0,\alpha_+\}$, the usual leverage effect
holds for large shocks (at least larger than 1) but is reversed for small ones. A measure of
the average leverage effect can be defined through the covariance
between $\eta_{t-1}$ and the current log-volatility.
We restrict our study to the case $p=q=1$, omitting subscripts to
simplify notation.
\begin{prop}[{\bf Leverage effect in the log-GARCH(1,1) model}]\label{prlev}
Consider the log-GARCH(1,1) model under the condition
$\rho(\bA^{(\infty)})<1$. Assume that the innovations $\eta_t$ are
symmetrically distributed, $E[|\log \eta_0|^{2}]<\infty$ and
$|\beta|+\frac12 (|\alpha_{+}|+|\alpha_{-}|)<1$. Then
\begin{equation}\label{leverage}
\mbox{cov}(\eta_{t-1},\log\sigma_t^{2})=\frac12
(\alpha_{+}-\alpha_{-})\left\{ E(|\eta_0|)
\tau+E(|\eta_0|\log\eta^2_{0}) \right\},\end{equation} where
$$\tau=E\log \sigma_t^2=
\frac{\omega + \frac12 (\alpha_{+}-\alpha_{-})E(\log\eta^2_{0})}{1-\beta-\frac12
(\alpha_{+}+\alpha_{-})}.
$$
\end{prop}
Thus, if the left hand side of (\ref{leverage}) is negative the
leverage effect is present: past negative innovations tend to
increase the log-volatility, and hence the volatility, more than
past positive innovations. However, the sign of the covariance is
more complicated to determine than for other asymmetric models: it
depends on all the GARCH coefficients, but also on the properties of
the innovations distribution. Interestingly, the leverage effect may
hold with $\alpha_{+}>\alpha_{-}.$ 

\section{Tail properties of the log-GARCH}
\label{secdep}

In this section, we investigate differences between the EGARCH and
the log-GARCH in terms of tail properties.

\subsection{Existence of moments}
We start by characterizing the existence of  moments for the
log-GARCH. The following result is an extension of Theorem 1 in
Bauwens {\it et al.}, 2008, to the asymmetric case (see also Theorem
2 in He {\it et al.}, 2002 for the symmetric case with $p=q=1$).
\begin{prop}[{\bf Existence of moments}]\label{prmom}
Assume that $\gamma ({\bf C})<0$ and  that
$\rho\left(\bA^{(\infty)}\right)<1$. Letting $\lambda=\max_{1\le
i\le q}\{|\alpha_{i+}|\vee |\alpha_{i-}|\}\sum_{\ell\ge
0}\|(\bA^{(\infty)})^{\ell}\|<\infty$, assume that for some $s>0$
\begin{equation}\label{condmom}
E\Big[\exp\Big\{ s \Big(\lambda  \vee 1\Big)|\log\eta_0^2|\Big\}\Big]<\infty,
\end{equation}
then the solution of the log-GARCH(p,q) model satisfies
$E|\epsilon_0|^{2s}<\infty$.
\end{prop}
\begin{rem}\label{remmom} {\rm In the case $p=q=1$, condition \eqref{condmom}
has a simpler form:
$$
E\Big[\exp\Big\{ s \Big(\frac{|\alpha_{1+}|\vee |\alpha_{1-}|}{1-|\alpha_{1+}+\beta_1|\vee |\alpha_{1-}+\beta_1|}\vee 1\Big)|\log\eta_0^2|\Big\}\Big]<\infty.
$$
} \end{rem} The following result provides a sufficient condition for
the Cramer's type condition \eqref{condmom}.
\begin{prop}\label{propcc}
If $E (|\eta_0|^s)<\infty$ for some $s>0$ and $\eta_0$ admits a
density $f$ around $0$  such that $f(y^{-1})=o(|y|^\delta)$ for
$\delta<1$ when $|y|\to\infty$ then
$E\exp(s_1|\log\eta_0^2|)<\infty$ for some $s_1>0$.
\end{prop}
In the case $p=q=1$, a very simple moment condition is given by the
following result.
\begin{prop}[{\bf Moment condition for the log-GARCH(1,1) model}]\label{prrvbis}
 Consider the log-GARCH(1,1) model.
 Assume that $E\log^+|\log\eta_0^2|<\infty$ and $E|\eta_0|^{2s}<\infty$ 
for $s>0$. Assume 
$\beta_1+\alpha_{1+}\in (0,1), \beta_1+\alpha_{1-}\in (0,1)$ and
$\alpha_{1+}\wedge \alpha_{1-}>0.$
 Then $\sigma^2_0$ and $\epsilon_0$ have
finite moments of order  $s/(\alpha_{1+}\vee\alpha_{1-})$ and
$2s/(\alpha_{1+}\vee\alpha_{1-}\vee 1)$ respectively.
\end{prop}
It can be noted that, for a given log-GARCH(1,1) process, moments
may exist at an arbitrarily large order. In this respect, log-GARCH
differ from standard GARCH and other GARCH specifications. In such
models the region of the parameter space such that $m$-th order
moment exist reduces to the empty set as $m$ increases.  For an
explicit expression of the unconditional moments in the case of
symmetric log-GARCH($p,q$) models,
 we refer the reader to Bauwens {\it et al.} (2008).\\


\subsection{Regular variation of the log-GARCH(1,1)} Under the
assumptions of Proposition \ref{prlev} we have an explicit
expression of the stationary solution. Thus it is possible to
establish the regular variation properties of the log-GARCH model.
Recall that $L$ is a slowly varying function iff $L(xy)/L(x)\to 1$
as $x\to \infty$ for any $y>0$. A random variable $X$ is said to be
regularly varying of index $s>0$ if there exists a slowly varying
function $L$ and $\tau\in [0,1]$ such that
$$
P(X>x)\sim \tau x^{-s}L(x)\quad\mbox{and}\quad P(X\le -x)\sim (1-\tau)x^{-s}L(x)\quad x\to +\infty.
$$
The following proposition asserts the regular  variation properties of the stationary solution of the log-GARCH(1,1) model.
\begin{prop}[Regular variation of the log-GARCH(1,1) model]\label{prrv}
Consider the log-GARCH(1,1) model.
 Assume that $E\log^+|\log\eta_0^2|<\infty$
that $\eta_0 $ is regularly varying with index $2s'>0$.
Assume 
$\beta_1+\alpha_{1+}\in (0,1), \beta_1+\alpha_{1-}\in (0,1)$ and
$\alpha_{1+}\wedge \alpha_{1-}>0.$ Then $\sigma^2_0$ and
$\epsilon_0$ are regularly varying with index
$s'/(\alpha_{1+}\vee\alpha_{1-})$ and
$2s'/(\alpha_{1+}\vee\alpha_{1-}\vee 1)$ respectively.
\end{prop}
The square root of the volatility, $\sigma_0$, thus have heavier
tails than the innovations when $\alpha_{1+}\vee\alpha_{1-}>1$.
Similarly, in  the EGARCH(1,1) model the observations can have a
much heavier tail than the innovations.
Moreover, when the innovations are light tailed distributed (for instance exponentially distributed),
the EGARCH can exhibit regular variation properties.
It is not the case for the log-GARCH(1,1) model.\\ 

In this context of heavy tail, a natural way to deal with the
dependence structure is to study the multivariate regular variation
of a trajectory. As the innovations are independent, the dependence
structure can only come from the volatility process. However, it is
also independent in the extremes. The following is a straightforward
application of Lemma 3.4 of Mikosch and Rezapur (2012).
\begin{prop}[Multivariate regular variation of the log-GARCH(1,1) model]
Assume the conditions of Proposition \ref{prrv} satisfied. Then the
sequence $(\sigma_t^2)$ is regularly varying with index
$s'/(\alpha_{1+}\vee\alpha_{1-})$. The limit measure of the vector
${\Sigma}_d^2= (\sigma_1^2,\ldots,\sigma_d^2)'$ is given by the
following limiting relation on the Borel $\sigma$-field of
$(R\cup\{+\infty\})^d\slash\{0_d\}$
$$
\frac{P(x^{-1}{\Sigma}_d^2\in\cdot)}{
P(\sigma^2> x)}\to
\frac{s'}{\alpha_{1+}\vee\alpha_{1-}}
\sum_{i=1}^d
\int_1^\infty y^{-s'/(\alpha_{1+}\vee\alpha_{1-})-1}1_{\{ye_i\in\cdot \}}dy,\quad x\to\infty.
$$
where $e_i$ is the $i$-th unit vector in $R^d$ and the convergence holds vaguely.
\end{prop}
As for the innovations, the  limiting measure above is concentrated
on the axes. Thus it is also the case for the log-GARCH(1,1) process
and its extremes values do not cluster. It is a drawback for
modeling stock returns when clusters of volatilities are stylized
facts. This lack of clustering is also observed for the EGARCH(1,1)
model in Mikosch and Rezapur (2012), in contrast with the GARCH(1,1)
model,
see Mikosch and Starica (2000).\\

\section{Estimating the log-GARCH by QML}
\label{sec3}

We now consider the statistical inference. Let $\epsilon_1,\dots,\epsilon_n$ be observations of the stationary solution of (\ref{logGARCH}),
where $\btheta$ is equal to an unknown value $\btheta_0$ belonging to some parameter space $\Theta\subset \mathbb{R}^d$, with $d=2q+p+1$.  A QMLE of $\btheta_0$ is defined as any measurable solution
$\widehat{\btheta}_n$  of
\begin{equation}\label{qml}
  \widehat{\btheta}_n=\argmin_{\btheta\in\Theta}\widetilde{Q}_n(\btheta),
\end{equation}
with
$$
  \widetilde{Q}_n(\btheta) =n^{-1}\sum_{t=r_0+1}^n\widetilde{\ell}_t(\btheta),\qquad
 \widetilde{\ell}_t(\btheta)= \frac{\epsilon_t^2}{\widetilde{\sigma}_t^2(\btheta)} +\log \widetilde{\sigma}_t^2(\btheta),
$$
where $r_0$ is a fixed integer and $\log \widetilde{\sigma}_t^2(\btheta)$ is recursively defined, for $t=1,2,\dots,n$, by
$$\log \widetilde{\sigma}_t^2(\btheta)=\omega+\sum_{i=1}^q\left(\alpha_{i+}1_{\{\epsilon_{t-i}>0\}}+\alpha_{i-}1_{\{\epsilon_{t-i}<0\}}\right)\log\epsilon_{t-i}^2
+\sum_{j=1}^p\beta_j\log \widetilde{\sigma}_{t-j}^2(\btheta),$$
using positive initial values for $\epsilon_0^2,\dots, \epsilon_{1-q}^2,\widetilde{\sigma}_{0}^2(\btheta),\dots, ,\widetilde{\sigma}_{1-p}^2(\btheta)$.
\begin{rem}[{\bf On the choice of the initial values}]
\label{valinit0} {\rm
It will be shown in the sequel that the choice of $r_0$ and of the
initial values is unimportant for the asymptotic behavior of the
QMLE, provided $r_0$ is fixed and there exists a real random
variable $K$ independent of $n$ such that
\begin{equation}
\label{condvi}
\sup_{\btheta\in\Theta}\left|\log\sigma_t^2(\btheta)-\log\widetilde{\sigma}_t^2(\btheta)\right|<K,\quad\mbox{a.s. for }t=q-p+1,\dots, q,
\end{equation}
where $\sigma_t^2(\btheta)$ is defined by (\ref{sigmadef}) below.
These conditions are supposed to hold in the sequel.}
\end{rem}
\begin{rem}[{\bf The empirical treatment of null returns}]
\label{nulreturn} {\rm Under the assumptions of
Theorem~\ref{theostatio}, almost surely  $\epsilon_t^2\neq 0.$
However, it may happen that some observations are equal to zero  or
are so close to zero that $\widehat{\btheta}_n$ cannot be computed
(the computation of the $\log\epsilon_t^2$'s being required). To
solve this potential problem, we imposed a lower bound for the
$|\epsilon_t|$'s. We took the lower bound $10^{-8}$, which is well
inferior to a beep point, and we checked that nothing was changed in
the numerical illustrations presented here when this lower bound was
multiplied or divided by a factor of 100.}
\end{rem}
We now need to introduce some notation. For any $\btheta\in\Theta$, let the polynomials ${\cal A}^+_{\btheta}(z)=\sum_{i=1}^q\alpha_{i,+}z^{i}$,
${\cal A}^-_{\btheta}(z)=\sum_{i=1}^q\alpha_{i,-}z^{i}$ and
${\cal B}_{\btheta}(z)=1-\sum_{j=1}^p\beta_{j}z^{j}.$ By
convention, ${\cal A}^{+}_{\btheta}(z)=0$ and ${\cal A}^{-}_{\btheta}(z)=0$ if $q=0$, and  ${\cal
B}_{\btheta}(z)=1$ if $p=0$. We also write ${\bf C}(\btheta_0)$ instead of ${\bf C}$ to emphasize  that the unknown parameter is $\btheta_0$.
The following assumptions are used to show the strong consistency of the QMLE.
\begin{itemize}
\item[\hspace*{1em} {\bf A1:}]
\hspace*{1em} $ \btheta_0\in\Theta$  and  $\Theta$ is compact.
\item[\hspace*{1em} {\bf A2:}]
\hspace*{1em} $\gamma \left\{{\bf C(\btheta_0)}\right\}<0 \quad
$ and $\quad \forall \btheta\in\Theta,\quad |{\cal
B}_{\btheta}(z)|=0 \Rightarrow |z|>1.$
\item[\hspace*{1em} {\bf A3:}]
\hspace*{1em} 
The support of $\eta_0$ contains at least two positive values
and two negative values, $E\eta_0^2=1$ and
$E|\log\eta_0^2|^{s_0}<\infty$ for some $s_0>0$.
\item[\hspace*{1em} {\bf A4:}] \hspace*{1em} If $p>0$,
${\cal A}^+_{\btheta_0}(z)$ and ${\cal A}^-_{\btheta_0}(z)$  have no common root with  ${\cal B}_{\btheta_0}(z)$. Moreover ${\cal A}^+_{\btheta_0}(1)+{\cal A}^-_{\btheta_0}(1)\neq
0$ and  $|\alpha_{0q+}|+|\alpha_{0q+}|+|\beta_{0p}|\neq 0$.
\item[\hspace*{1em} {\bf A5:}]
\hspace*{1em} $E\left|\log \epsilon_t^2\right|<\infty$.
\end{itemize}
Assumptions {\bf A1, A2} and {\bf A4} are similar to those required
for the consistency of the QMLE in standard GARCH models (see Berkes
et al. 2003, Francq and Zakoian, 2004). Assumption {\bf A3}
precludes a mass at zero for the innovation, and, for
identifiability reasons, imposes non degeneracy of the positive and
negative parts of $\eta_0$. Note that, for other GARCH-type models,
the absence of a lower bound for the volatility can entail
inconsistency of the (Q)MLE (see Francq and Zakoïan (2010b) for a
study of a finite-order version of the LARCH($\infty$) model
introduced by Robinson (1991)). This is not the case for the
log-GARCH under {\bf A5}.
 Note that this assumption can be
replaced by the sufficient conditions given in
Proposition~\ref{logmoment} (see also
Example~\ref{examplelog-garch}).
\begin{theo}[{\bf Strong consistency of the QMLE}]\label{consistency}
Let $(\widehat{\btheta}_n)$ be a sequence of QMLE satisfying (\ref{qml}), where the $\epsilon_t$'s follow the asymmetric log-GARCH model of parameter $\btheta_0$.
Under the assumptions (\ref{condvi}) and {\bf A1}-{\bf A5}, almost surely $\widehat{\btheta}_n\to \btheta_0$ as $n\to\infty$.
\end{theo}
Let us now study the asymptotic normality of the QMLE. We need the
classical additional assumption:
\begin{itemize}
\item[\hspace*{1em} {\bf A6:}]
\hspace*{1em} $\btheta_0\in\stackrel{\circ}{\Theta}$ and
$\kappa_4:=E(\eta_0^4)<\infty$.
\end{itemize}
Because the volatility $\sigma_t^2$ is not bounded away from $0$, we also need the following non classical assumption.
\begin{itemize}
\item[\hspace*{1em} {\bf A7:}] There exists $s_1>0$ such that
$E\exp(s_1|\log\eta_0^2|)<\infty$, and $\rho(\bA^{(\infty)})<1$.
\end{itemize}
The Cramer condition on $|\log \eta_0^2|$ in {\bf A7} is verified if
$\eta_t$ admits a density $f$ around $0$ that does not explode too
fast (see Proposition \ref{propcc}).

Let $\nabla Q=(\nabla_1 Q,\dots, \nabla_d Q)'$ and $\mathbb H
Q=(\mathbb H_{1.} Q',\dots, \mathbb H_{d.}Q')'$ be the vector and
matrix of the first-order and second-order partial derivatives of a
function $Q:\Theta\to \mathbb{R}$.
\begin{theo}[{\bf Asymptotic normality of the QMLE}]\label{normality}
Under the assumptions of Theorem~\ref{consistency}
 and {\bf A6}-{\bf A7}, we have $\sqrt n (\widehat\btheta_n-\btheta_0)\stackrel{d}{\to}\mathcal N(\bzero,(\kappa_4-1)\bf J^{-1})$ as $n\to\infty$, where
 ${\bf J}=E[\nabla \log \sigma_t^2(\btheta_0) \nabla \log
\sigma_t^2(\btheta_0)']$ is a positive definite matrix and
 $\stackrel{d}{\to}$ denotes convergence in distribution.
\end{theo}
It is worth noting that for the general EGARCH model, no similar
results, establishing the consistency and the asymptotic normality,
exist. See however Wintenberger (2013) for the EGARCH(1,1). The
difficulty with the EGARCH is to  invert the volatility, that is to
write $\sigma_t^2(\btheta)$ as a well-defined function of the past
observables. In the log-GARCH model,  invertibility reduces to the
standard assumption on ${\cal B}_{\btheta}$ given in {\bf A2}.

\section{Asymmetric log-ACD model for duration data}

The dynamics of duration between stock price changes has attracted
much attention in the econometrics literature. Engle and Russel
(1997) proposed the Autoregressive Conditional Duration (ACD) model,
which assumes that the duration between price changes has the
dynamics of the square of a GARCH. Bauwens and Giot (2000 and 2003)
introduced logarithmic versions of the ACD, that do not constrain
the sign of the coefficients (see also Bauwens, Giot, Grammig and
Veredas (2004) and Allen, Chan, McAleer and Peiris (2008)). The
asymmetric ACD of Bauwens and Giot (2003) applies to pairs of
observation $(x_i,y_i)$, where $x_i$ is the duration between two
changes of the bid-ask quotes posted by a market maker and $y_i$ is
a variable indicating the direction of change of the mid price
defined as the average of the bid and ask prices ($y_i=1$ if the mid
price increased over duration $x_i$, and $y_i=-1$ otherwise). The
asymmetric log-ACD proposed by Bauwens and Giot (2003) can be
written as
\begin{equation}\label{logacd}
\left\{\begin{array}{lll}  x_i&=&\psi_i z_i, \\
\log \psi_i&=&\omega+ \sum_{k=1}^q \left(\alpha_{k+}1_{\{y_{i-k}=1\}}+\alpha_{k-}1_{\{y_{i-k}=-1\}}\right)\log x_{i-k}\\
&&+ \sum_{j=1}^p \beta_{j}\log \psi_{i-j},\end{array}\right.
\end{equation}
where $(z_i)$ is an iid sequence of positive variables with mean 1 (so that $\psi_i$ can be interpreted as the conditional mean of the duration $x_i$).
Note that $\epsilon_t:=\sqrt{x_t}y_t$ follows the log-GARCH model (\ref{logGARCH}), with $\eta_t=\sqrt{z_t}y_t$. Consequently, the results of the present paper also apply to log-ACD models. In particular, the parameters of (\ref{logacd}) can be estimated by fitting model (\ref{logGARCH})
on $\epsilon_t=\sqrt{x_t}y_t$.

\section{Numerical Applications}
\label{sec7}

\subsection{An application to exchange rates} We consider
returns series of the daily exchange rates of the American Dollar
(USD), the Japanese Yen (JPY), the British Pound (BGP), the Swiss
Franc (CHF)  and Canadian Dollar (CAD) with respect to the Euro. The
observations cover the period from January 5, 1999 to January 18,
2012, which corresponds to 3344 observations. The data were obtained
from the web site

\noindent\url{http://www.ecb.int/stats/exchange/eurofxref/html/index.en.html}.

Table~\ref{TauxLGARCH} displays the estimated log-GARCH(1,1) and
EGARCH(1,1) models for each series. As in Wintenberger (2013), the
optimization of the EGARCH(1,1) models has been performed under the
constraints $\delta\geq |\gamma|$ and
$$\sum_{t=1}^n\log\left[\max\left\{\beta,\frac{1}{2}\left(\gamma\epsilon_{t-1}+\delta|\epsilon_{t-1}|\right)\exp\left(-\frac{1}{2}\frac{\alpha}{1-\beta}\right)\right\}-\beta\right]<0.$$
These constraints guarantee the invertibility of the model, which is
essential to obtain a model that can be safely used for prediction
(see Wintenberger  (2013) for details).
For all series, except the CHF,
condition (\ref{allmombis}) ensuring the existence of any log-moment for the log-GARCH
is satisfied.
For all models, the persistence parameter $\beta$ is very high. The
last column  shows that for the USD and the GBP, the log-GARCH has a
higher (quasi) log-likelihood than the EGARCH. The converse is true
for the three other assets. A study of the residuals, not reported
here, is in accordance with the better fit of one particular model
for each series. It is also interesting to see that the two models
do not detect asymmetry for the same series. Moreover, models for
which the symmetry assumption is rejected (EGARCH for the JPY and
CHF, log-GARCH for the USD series) is also the preferred one
in terms of log-likelihood. 
This study
confirms that the models do not capture exactly the same empirical
properties, and are thus not perfectly substitutable.

\renewcommand{\baselinestretch}{1}

\begin{table}
\begin{center}\caption{\label{TauxLGARCH}
Log-GARCH(1,1) and EGARCH(1,1) models fitted by QMLE on daily
returns of exchange rates. The estimated standard deviation are
displayed into
brackets. The 6th column gives the $p$-values of the Wald test for symmetry ($\alpha_+=\alpha_-$ for the log-GARCH and $\gamma=0$ for the  EGARCH),
in bold face when the null hypothesis is rejected at level greater than 1\%. The last column gives the log-likelihoods (up to a constant) for the two models
with the largest in bold face.}
\begin{tabular}{lcccccc}
\vspace*{0.0cm}\\
\hline\hline\vspace*{-0.3cm}\\
\multicolumn{7}{l}{Log-GARCH} \\
 &  $\widehat{\omega}$ & $\widehat{\alpha}_+$ & $\widehat{\alpha}_-$ & $\widehat{\beta}$&p-val& Log-Lik.  \\
 USD  &    0.024     (0.005)    &    0.027       (0.004)     &     0.016      (0.004)     &     0.971      (0.005)& {\bf 0.01}  &   \textbf{-0.104}\\
 JPY  &    0.051     (0.007)    &    0.037       (0.006)     &     0.042      (0.006)     &     0.952      (0.006)& 0.36  &   -0.354         \\
 GBP  &    0.032     (0.006)    &    0.030       (0.005)     &     0.029      (0.005)     &     0.964      (0.006)& 0.84  &   \textbf{0.547} \\
 CHF  &    0.057     (0.012)    &    0.046       (0.008)     &     0.036      (0.007)     &     0.954      (0.008)& 0.11  &   1.477          \\
 CAD  &    0.021     (0.005)    &    0.025       (0.004)     &     0.017      (0.004)     &     0.969      (0.006)& 0.12   &  -0.170         \\
\vspace*{0.1cm}\\
\multicolumn{7}{l}{EGARCH} \\
 &  $\widehat{\omega}$ & $\widehat{\gamma}$ & $\widehat{\delta}$ & $\widehat{\beta}$&p-val& Log-Lik. \\
USD  &    -0.172      (0.027)    &    -0.014      (0.013)     &     0.189      (0.029)     &     0.970      (0.009)& 0.28 &   -0.110           \\
JPY  &    -0.209      (0.025)    &    -0.091      (0.016)     &     0.236      (0.029)     &     0.955      (0.008)&  {\bf 0.00} &  \textbf{-0.342}   \\
 GBP  &    -0.242     (0.035)    &    -0.019      (0.016)     &     0.233      (0.032)     &     0.959      (0.009)&  0.24 &   0.540            \\
 CHF  &    -0.103     (0.021)    &    -0.046      (0.014)     &     0.087      (0.018)     &     0.986      (0.004)&  {\bf 0.00} &  \textbf{1.575}    \\
 CAD  &    -0.067     (0.013)    &    -0.005      (0.009)     &     0.076      (0.015)     &     0.990      (0.004)&  0.57 &  \textbf{-0.160}   \\
\hline
\end{tabular}
\end{center}
\end{table}
\renewcommand{\baselinestretch}{1.5}

\subsection{A Monte Carlo experiment}
To evaluate the finite sample performance of the QML for the two
models we made the following numerical experiments. We first
simulated the log-GARCH(1,1) model, with $n=3344$, $\eta_t\sim {\cal
N}(0,1)$, and a parameter close to those 
of Table~\ref{TauxLGARCH}, that is
$\btheta_0=(0.024,0.027,0.016,0.971)$. Notice that assumptions {\bf
A1}--{\bf A4} required for the consistency are clearly satisfied.
Since
$|\beta_0|+\frac{1}{2}\left(|\alpha_{0+}|+|\alpha_{0-}|\right)<1$,
{\bf A5} is also satisfied in view of
Example~\ref{examplelog-garch}. The assumptions {\bf A6}-{\bf A7}
required for the asymptotic normality are also satisfied, noting
that $|\alpha_{0+}+\beta_0|\vee|\alpha_{0-}+\beta_0|<1$ and using
Proposition~\ref{anylogmoment}. The first part of
Table~\ref{asTauxLGARCH} displays the log-GARCH(1,1) models fitted
on these simulations. This table shows that the log-GARCH(1,1) is
accurately estimated. Note that the estimated models satisfy also
the assumptions {\bf A1}-{\bf A7} used to show the consistency and
asymptotic normality. We also estimated EGARCH(1,1) models on the
same simulations. The results are presented in the second part of
Table~\ref{asTauxLGARCH}.  Comparing the log-likelihood given in the
last column of Table~\ref{asTauxLGARCH}, one can see that, as
expected,  the likelihood of the log-GARCH model is greater than
that of the (misspecified) EGARCH model, for all the simulations.

In a second time, we repeated the same experiments for simulations
of an EGARCH(1,1) model of parameter $(\omega_0, \gamma_0, \delta_0,
\beta_0)=(-0.204,-0.012,0.227,0.963)$. Table~\ref{asTauxLGARCH2} is
the analog of   Table~\ref{asTauxLGARCH} for the simulations of this
EGARCH model instead of the log-GARCH. The EGARCH are satisfactorily
estimated, and, once again, the simulated model has a higher
likelihood than the misspecified model. From this simulation
experiment, we draw the conclusion that it makes sense to select the
model with the higher likelihood, as we did for the series of
exchange rates.

\renewcommand{\baselinestretch}{1}

\begin{table}
\begin{center}\caption{\label{asTauxLGARCH}
Log-GARCH(1,1) and EGARCH(1,1)  models fitted on 5 simulations of a
log-GARCH(1,1) model. The estimated standard deviation are displayed
into brackets. The  larger log-likelihood is displayed in bold
face.}
\begin{tabular}{lccccc}
\vspace*{0.0cm}\\
\hline\hline\vspace*{-0.3cm}\\
\multicolumn{6}{l}{Log-GARCH} \\
Iter &  $\widehat{\omega}$ & $\widehat{\alpha}_+$ & $\widehat{\alpha}_-$ & $\widehat{\beta}$& Log-Lik.\\
 1  &  0.025   (0.004)    &    0.028    (0.004)     &     0.018    (0.004)     &     0.968     (0.005)    &   \textbf{-0.415}    \\
 2  &  0.021   (0.003)    &    0.023    (0.003)     &     0.013    (0.003)     &     0.976     (0.004)    &   \textbf{-0.634}    \\
 3  &  0.026   (0.003)    &    0.028    (0.004)     &     0.017    (0.003)     &     0.969     (0.004)    &   \textbf{-0.754}    \\
 4  &  0.022   (0.003)    &    0.024    (0.004)     &     0.018    (0.003)     &     0.972     (0.004)    &   \textbf{-0.389}    \\
 5  &  0.024   (0.003)    &    0.028    (0.004)     &     0.014    (0.003)     &     0.974     (0.003)    &   \textbf{-0.822}    \\
\vspace*{0.1cm}\\
\multicolumn{6}{l}{EGARCH} \\
Iter  &  $\widehat{\omega}$ & $\widehat{\gamma}$ & $\widehat{\delta}$ & $\widehat{\beta}$& Log-Lik.\\
 1  &   -0.095    (0.016)    &    -0.014     (0.009)                &     0.104     (0.017)     &     0.976      (0.006)   & -0.424        \\
 2  &   -0.127    (0.018)    &    \phantom{-}0.009      (0.010)     &     0.148     (0.021)     &     0.976      (0.007)   & -0.645        \\
 3  &   -0.147    (0.018)    &    \phantom{-}0.001      (0.010)     &     0.177     (0.022)     &     0.971      (0.007)   & -0.770        \\
 4  &   -0.136    (0.019)    &    -0.012     (0.010)                &     0.155     (0.022)     &     0.976      (0.007)   & -0.404        \\
 5  &   -0.146    (0.019)    &    -0.009     (0.010)                &     0.177     (0.023)     &     0.971      (0.007)   & -0.842        \\
\hline
\end{tabular}
\end{center}
\end{table}

%
\begin{table}
\begin{center}\caption{\label{asTauxLGARCH2}
As Table~\ref{asTauxLGARCH}, but for 5 simulations of an EGARCH(1,1)
model.}
\begin{tabular}{lccccc}
\vspace*{0.0cm}\\
\hline\hline\vspace*{-0.3cm}\\
\multicolumn{6}{l}{Log-GARCH} \\
Iter &  $\widehat{\omega}$ & $\widehat{\alpha}_+$ & $\widehat{\alpha}_-$ & $\widehat{\beta}$& Log-Lik.\\
 1  &   0.039   (0.008)   &    0.071    (0.008)     &     0.052     (0.007)     &     0.874    (0.015) &   -0.350    \\
 2  &   0.055   (0.006)   &    0.058    (0.007)     &     0.052     (0.006)     &     0.913    (0.010) &   -0.476    \\
 3  &   0.052   (0.008)   &    0.070    (0.008)     &     0.060     (0.007)     &     0.873    (0.015) &   -0.468    \\
 4  &   0.051   (0.008)   &    0.076    (0.008)     &     0.056     (0.007)     &     0.878    (0.014) &   -0.416    \\
 5  &   0.056   (0.007)   &    0.061    (0.007)     &     0.060     (0.007)     &     0.896    (0.012) &   -0.517    \\
\vspace*{0.1cm}\\
\multicolumn{6}{l}{EGARCH} \\
Iter  &  $\widehat{\omega}$ & $\widehat{\gamma}$ & $\widehat{\delta}$ & $\widehat{\beta}$& Log-Lik.\\
 1  & -0.220   (0.022)    &    -0.024   (0.013)     &     0.235   (0.023)     &     0.950   (0.010)  &-\textbf{0.335}    \\
 2  & -0.196   (0.020)    &    -0.029   (0.012)     &     0.219   (0.022)     &     0.961   (0.008)  &\textbf{-0.468}    \\
 3  & -0.222   (0.022)    &    -0.005   (0.013)     &     0.241   (0.024)     &     0.947   (0.010)  &\textbf{-0.448}    \\
 4  & -0.227   (0.022)    &    -0.025   (0.012)     &     0.248   (0.023)     &     0.950   (0.010)  &\textbf{-0.402}    \\
 5  & -0.209   (0.021)    &    -0.003   (0.012)     &     0.234   (0.023)     &     0.955   (0.009)  &\textbf{-0.504}    \\
 \hline
\end{tabular}
\end{center}
\end{table}
\renewcommand{\baselinestretch}{1.5}


\section{Proofs}
\label{secProofs}

\subsection{Proof of Theorem \ref{theostatio}}
Since the random variable $\|\bC_0\|$ is bounded, we have
$E\log^+\|\bC_0\|<\infty$. The moment condition on $\eta_t$ entails
that we also have $E\log^+\|\bb_0\|<\infty$. When $\gamma ({\bf C})
<0$, Cauchy's root test shows that, almost surely (a.s.), the series
\begin{equation}\label{vecinfmul}
\bz_t=\bb_t+\Sum_{n=0}^{\infty}
\bC_t\bC_{t-1}\cdots \bC_{t-n} \bb_{t-n-1}\end{equation}
converges absolutely for all $t$ and satisfies (\ref{vectorform}). A
strictly stationary solution to model (\ref{logGARCH}) is then
obtained as $\epsilon_t=
\exp\left\{\frac{1}{2}\bz_{2q+1,t}\right\}\eta_t$, where $\bz_{i,t}$
denotes the $i$-th element of $\bz_t$. This solution is non
anticipative and  ergodic, as a measurable function of
$\{\eta_u,u\leq t\}$.

We now prove that (\ref{vecinfmul}) is the  unique nonanticipative
solution of (\ref{vectorform}) when $\gamma ({\bf C})< 0$.   Let
$(\bz^*_t)$ be a strictly stationary process satisfying
$\bz^*_t=\bC_t\bz^*_{t-1}+\bb_t$. For all $N\geq 0$,
$$\bz^*_t=\bz_t(N)+
\bC_t\ldots \bC_{t-N}\bz^*_{t-N-1}, \quad
\bz_t(N)=\bb_t+\Sum_{n=0}^{N}
\bC_t\bC_{t-1}\cdots \bC_{t-n} \bb_{t-n-1}.
$$
We then have
$$\|\bz_t-\bz^*_t\|\leq \left\|\Sum_{n=N+1}^{\infty}
\bC_t\bC_{t-1}\cdots \bC_{t-n} \bb_{t-n-1}\right\|+
\|\bC_t\ldots \bC_{t-N}\|\|\bz^*_{t-N-1}\|.$$ The first term
in the right-hand side tends to 0 a.s. when $N\to \infty.$ The
second term tends to 0 in probability because  $\gamma ({\bf C}) <0$
entails that $\|\bC_t\ldots \bC_{t-N}\|\to 0$  a.s.  and the
distribution of $\|\bz^*_{t-N-1}\|$ is independent of $N$ by
stationarity.  We have shown that $\bz_t-\bz^*_t\to 0$ in
probability when $N\to \infty.$ This quantity being independent of
$N$ we have $\bz_t=\bz^*_t$ a.s. for any $t$. \zak

\subsection{Proof of Proposition~\ref{fracmoment}} Let $X$ be a
random variable such that $X>0$ a.s. and $EX^r<\infty$ for some
$r>0$. If $E\log X<0$, then there exists $s>0$ such that $EX^s<1$
(see {\it  e.g.} Lemma 2.3 in Berkes, Horv\'ath and  Kokoszka,
2003). Noting that $E\left\|\bC_t\cdots\bC_1\right\|\leq
(E\left\|\bC_1\right\|)^t<\infty$ for all $t$, the previous result
shows that when $\gamma ({\bf C})<0$ we have
$E\|\bC_{k_0}\cdots\bC_1\|^s<1$ for some $s>0$ and some $k_0\geq 1$.
One can always assume that $s<1$. In view of (\ref{vecinfmul}), the
$c_r$-inequality and standard arguments (see {\it  e.g.} Corollary
2.3 in Francq and Zakoïan, 2010a) entail that $E\|\bz_t\|^s<\infty$,
provided $E\|\bb_t\|^s<\infty$,  which holds true when $s\leq s_0$.
The conclusion follows.\zak

\subsection{Proof of Proposition~\ref{logmoment}} By
(\ref{Markovian}),  componentwise we have
\begin{equation}\label{eq:abs}\mathrm{Abs}(\bsigma_{t,r})\leq \mathrm{Abs}({\bu}_t)+\sum_{\ell=0}^{\infty}\bA_{t,\ell}\mathrm{Abs}(\bu_{t-\ell-1}),
\quad \bA_{t,\ell}:=\prod_{j=0}^{\ell}\mathrm{Abs}({\bA}_{t-j}),\end{equation}
where each element of the series is defined a priori in
$[0,\infty]$. In view of the form (\ref{MatA}) of the matrices
$\bA_t$, each element of
$$\bA_{t,\ell}\mathrm{Abs}(\bu_{t-\ell-1})=|u_{t-\ell-1}|\prod_{j=0}^{\ell}\mathrm{Abs}({\bA}_{t-j})\be_1$$ is  a sum of products of the form $|u_{t-\ell-1}|\prod_{j=0}^k|\mu_{\ell_j}(\eta_{t-i_j})|$
with $0\leq k\leq \ell$ and $0\leq i_0<\cdots<i_k\leq \ell+1$. To
give more detail, consider for instance the case $r=3$. We then have
$$\bA_{t,1}\mathrm{Abs}(\bu_{t-2})=
\left(\begin{array}{c}
|\mu_1(\eta_{t-1})||\mu_1(\eta_{t-2})||u_{t-2}|+|\mu_2(\eta_{t-2})||u_{t-2}|\\
|\mu_1(\eta_{t-2})||u_{t-2}|\\|u_{t-2}|
\end{array}\right).$$
Noting that $|u_{t-\ell-1}|$ is a function of $\eta_{t-\ell-2}$ and
its past values, we obtain
$E\bA_{t,1}\mathrm{Abs}(\bu_{t-2})=E\mathrm{Abs}({\bA}_{t})E\mathrm{Abs}({\bA}_{t-1})E\mathrm{Abs}(\bu_{t-2})$.
More generally, it can be shown by induction on $\ell$ that the
$i$-th element of the vector
$\bA_{t-1,\ell-1}\mathrm{Abs}(\bu_{t-\ell-1})$ is independent of the
$i$-th element of the first row of $\mathrm{Abs}(\bA_t)$. It follows
that
$E\bA_{t,\ell}\mathrm{Abs}(\bu_{t-\ell-1})=E\mathrm{Abs}(\bA_{t})E\bA_{t-1,\ell-1}\mathrm{Abs}(\bu_{t-\ell-1})$.
The property extends to $r\neq 3$. Therefore, although the matrices
involved in the product $\bA_{t,\ell}\mathrm{Abs}(\bu_{t-\ell-1})$
are not independent (in the case $r>1$), we have
\begin{eqnarray*}
E\bA_{t,\ell}\mathrm{Abs}(\bu_{t-\ell-1})
&=&\prod_{j=0}^{\ell}E\mathrm{Abs}({\bA}_{t-j})E\mathrm{Abs}(\bu_{t-\ell-1})
=\left(\bA^{(1)}\right)^{\ell+1}E\mathrm{Abs}(\bu_{1}).
\end{eqnarray*}
In view of (\ref{eq:abs}), the condition $\rho(\bA^{(1)})<1$ then
entails that $E\mathrm{Abs}(\bsigma_{t,r})$ is finite.

The case $r=1$ is treated by noting that
$\bA_{t,\ell}\mathrm{Abs}(\bu_{t-\ell-1})$ is a product of
independent  random variables.

To deal with the cases $r\neq 1$ and $m\neq 1$,  we work with
(\ref{vectorform}) instead of (\ref{Markovian}). This Markovian
representation has an higher dimension  but involves independent
coefficients $\bC_t$. Define $\bC_{t,\ell}$ by replacing $\bA_{t-j}$
by $\bC_{t-j}$ in $\bA_{t,\ell}$. We then have
$$E\bC_{t,\ell}^{\otimes m}\mathrm{Abs}(\bb_{t-\ell-1})^{\otimes m}=\left(\bC^{(m)}\right)^{\ell+1}E\mathrm{Abs}(\bb_{1})^{\otimes m}.$$
For all $m\geq 1$, let $\|\bM\|_m=(E\|\bM\|^m)^{1/m}$ where
$\|\bM\|$ is the sum of the absolute values of the elements of the
matrix $\bM$. Using the elementary relations
$\|\bM\|\|\bN\|=\|\bM\otimes\bN\|$ and
$E\|\mathrm{Abs}(\bM)\|=\|E\mathrm{Abs}(\bM)\|$ for any matrices
$\bM$ and $\bN$, the condition $\rho(\bC^{(m)})<1$  entails
$E\left\|\bC_{t,\ell}\mathrm{Abs}(\bb_{t-\ell-1})\right\|^m=\|E\bC_{t,\ell}^{\otimes
m}\mathrm{Abs}(\bb_{t-\ell-1})^{\otimes m}\|\to 0$ at the
exponential rate as $\ell\to\infty$, and thus
$$\left\|\mathrm{Abs}(\bz_t)\right\|_m\leq \left\|\mathrm{Abs}(\bb_t)\right\|_m+\sum_{\ell=0}^{\infty}\left\|\bC_{t,\ell}\mathrm{Abs}(\bb_{t-\ell-1})\right\|_m<\infty,$$
which allows to conclude.\zak

\subsection{Proof of
Proposition~\ref{anylogmoment}}
It follows from (\ref{eq:abs}) that componentwise we have
\begin{equation}\mathrm{Abs}(\bsigma_{t,r})\leq
\sum_{\ell=0}^{\infty}(\bA^{(\infty)})^{\ell}\mathrm{Abs}(\bu_{t-\ell}).\label{eq:abs1}
\end{equation} Therefore, the condition
$\rho(\bA^{(\infty)})<1$
ensures the existence of $E|\log\epsilon_t^2|^m$ at any order $m$,
provided $\gamma ({\bf C})<0$ and $E|\log\eta_0^2|^{m}<\infty$. Now
in view of the companion form of the matrix $\bA^{(\infty)}$ (see
{\it  e.g.} Corollary 2.2 in Francq and Zakoïan, 2010a), the
equivalence in (\ref{allmombis}) holds.\zak

\subsection{Proof of Proposition \ref{prlev}}

 By the concavity of the logarithm
function, the condition $|\alpha_++\beta||\alpha_-+\beta|<1$ is
satisfied.  By Example \ref{ex11} and the symmetry of the
distribution of $\eta_0$, the existence of a strictly stationary
solution
 $(\epsilon_t)$  satisfying $E|\log \epsilon_t^2 |<\infty$  is thus guaranteed.
Let
$$a_t=(\alpha_{+}1_{\{\eta_{t}>0\}}+\alpha_{-}1_{\{\eta_{t}<0\}})\eta_t, \qquad
b_t=(\alpha_{+}1_{\{\eta_{t}>0\}}+\alpha_{-}1_{\{\eta_{t}<0\}})\eta_t\log\eta^2_{t}.
$$
We have $Ea_t=(\alpha_{+}-\alpha_{-})E(\eta_01_{\{\eta_{0}>0\}})$
and
$Eb_t=(\alpha_{+}-\alpha_{-})E(\eta_0\log\eta^2_{0}1_{\{\eta_{0}>0\}}),
$ using the symmetry assumption for the second equality. Thus
\begin{eqnarray*} \mbox{cov}(\eta_{t-1},\log(\sigma_t^{2}))&=&
E[a_{t-1}\log(\sigma_{t-1}^{2})+b_{t-1}],
\end{eqnarray*}
and the conclusion follows. \zak

\subsection{Proof of Proposition \ref{prmom}}
By definition, $|\log(\sigma_t^2)|\le \|\bsigma_{t,r}\|=
\|\mathrm{Abs}(\bsigma_{t,r})\|$. Then,  we have
\begin{eqnarray*}
E|\sigma_t^2|^s&\leq & E\left\{\exp(s \|\mathrm{Abs}(\bsigma_{t,r})\|)\right\}= \sum_{k=0}^\infty \frac{s^k  \|\mathrm{Abs}(\bsigma_{t,r})\|_k^k}{k!}\\
&\leq&\sum_{k=0}^\infty \frac{s^k  \|\mathrm{Abs}({\bu}_{0})\|^k_k \left\{\sum_{\ell=0}^\infty    \|(\bA^{(\infty)})^{\ell}\|\right\}^k}{k!}\\
&=&E\exp\left\{s\|\mathrm{Abs}({\bu}_{0})\|\sum_{\ell=0}^\infty    \|(\bA^{(\infty)})^{\ell}\|\right\},
\end{eqnarray*}
where the last inequality comes from \eqref{eq:abs1}.
By definition ${\bu}_{0}=(u_0,0_{r-1}')'$ with
$$
u_0=\omega +\sum_{i=1}^q(\alpha_{i+}1_{\eta_{-i}>0}+\alpha_{i-}1_{\eta_{-i}<0})\log \eta_{-i}^2.
$$
Thus $\|\mathrm{Abs}({\bu}_{0})\|\le  |u_0|\leq |\omega|+\max_{1\leq
i\leq q}|\alpha_{i+}|\vee |\alpha_{i-}|\sum_{j=1}^q |\log
\eta_{-j}^2|$ and it follows that
\begin{eqnarray*}
E|\sigma_t^2|^s
&\leq&\exp\left\{s|\omega|\sum_{\ell=0}^\infty    \|(\bA^{(\infty)})^{\ell}\|\right\}\left\{E\exp\left(s\lambda |\log\eta^2_0|\right)\right\}^q<\infty
\end{eqnarray*}
under \eqref{condmom}. \zak

\subsection{Proof of Proposition \ref{propcc}} Without loss of
generality assume that $f$ exists on $[-1,1]$. Then there exists
$M>0$ such that $f(1/y)\le M|y|^\delta$ for all $y\ge 1$ and we
obtain
\begin{eqnarray*}
E\exp(s_1|\log\eta_0^2|)&\le& \int_{|x|<1}\exp(2s_1\log(1/x))f(x)dx+\int\exp(s_1\log(x^2))dP_\eta(x)\\
&\le& 2M\int_1^\infty y^{2(s_1-1)+\delta} dy+E(|\eta_0|^{2s_1}).
\end{eqnarray*}
The upper bound is finite for sufficiently small $s_1$ and the
result is proved. \zak

\subsection{Proof of Proposition \ref{prrvbis}}
We will use Tweedie's (1988) criterion, which we recall for the
reader's convenience. 
Let $(X_t)$ denote a temporally homogeneous Markov chain on a state
space $E$, endowed with a $\sigma$-field ${\cal F}$.
\begin{lem}[adapted from Tweedie (1988), Theorem 1]\label{Tweedielemme}
Suppose $\mu$ is a subinvariant measure, that is,
$$\mu(B)\geq \int_E \mu(dy)P(y,B), \quad \forall B\in {\cal F}$$
 and $A\in {\cal
F}$ is such that $0<\mu(A)<\infty.$ Suppose there exist a
nonnegative measurable function $V$ on $E$, and constants $K>0$ and
$c \in (0,1)$ such that
\begin{enumerate}[i)]
\item $E[V(X_t)\mid X_{t-1}=x]\leq K, \quad x\in A,$
\item $E[V(X_t)\mid X_{t-1}=x]\leq (1-c)V(x), \quad x\in A^c.$
\end{enumerate}
Then, $\mu$ is a finite invariant measure for the chain $(X_t)$ and
$\int_E Vd\mu<\infty.$
\end{lem}
\noindent It should be noted that this criterion does not make any
irreducibility assumption. We have
$u_t=\omega+(\alpha_{1+}1_{\{\eta_{t-1}>0\}}+\alpha_{1-}1_{\{\eta_{t-1}<0\}})\log(\eta^2_{t-1}),
$ and
$$\sigma_t^2=e^{u_t}\left(\sigma_{t-1}^{2(\beta+\alpha_{1+})}1_{\eta_{t-1}>0}+
\sigma_{t-1}^{2(\beta+\alpha_{1-})}1_{\eta_{t-1}<0}\right),$$ which
shows that $(\sigma_t^2)$ is a temporally homogeneous Markov chain
on $\mathbb{R}^{+*}$. 
By Example \ref{ex11}, the conditions ensuring the existence of a
strictly stationary solution are satisfied. The stationary
distribution thus defines an invariant probability $\mu$ for the
Markov chain $(\sigma_t^2)$. Let $A=[0,K]$ for some $K>0$ such that
$\mu(A)>0$.

The existence of a $s$-order moment for $\eta^2_{0}$ entails that
$e^{u_1}$ admits a moment of order
$s_0:=s/(\alpha_{1+}\vee\alpha_{1-})$.
Let $V(x)=1+x^{s_0}$ for $x>0.$ For any $x>0$ and for $0<c<1$, we
have for $x\in A^c$ with $K$ sufficiently large,
\begin{eqnarray*}
E[V(\sigma_t^2)\mid \sigma_{t-1}^2=x]&=&1+
x^{(\beta_1+\alpha_{1+})s_0}E[e^{s_0u_1}1_{\eta_0>0}]+
x^{(\beta_1+\alpha_{1-})s_0}E[e^{s_0u_1}1_{\eta_0<0}]\\
&\leq & (1-c)V(x),
\end{eqnarray*}
because $\beta_1+\alpha_{1+}<1$ and $\beta_1+\alpha_{1-}<1$. On the
other hand it is clear that $E[V(\sigma_t^2)\mid \sigma_{t-1}^2=x]$
is bounded for $x$ belonging to $A$. It follows by the above lemma
that $E_{\mu}[V(\sigma_t^2)]<\infty$ where the expectation is
computed with the stationary distribution.
 \zak

\subsection{Proof of Proposition \ref{prrv}}
To prove the first assertion, note that if $\eta_0$ is regularly
varying of index $2s'$ then $\eta_0^2$ is regularly varying of index
$s'$. Thus
$u_1=\omega+(\alpha_{1+}1_{\{\eta_{0}>0\}}+\alpha_{1-}1_{\{\eta_{0}<0\}})\log(\eta^2_{0})
$ is such that
\begin{eqnarray*}
P(e^{u_1}>x)&=& P(\eta_{0}>0)P\left\{(\eta_0^{2})^{\alpha_{1+}}>xe^{-\omega}\mid \eta_0>0\right\}\\
&&+P(\eta_{0}<0)P\left\{(\eta_0^{2})^{\alpha_{1-}}>xe^{-\omega}\mid \eta_0<0\right\}.
\end{eqnarray*}
Then $e^{u_1}$ is also regularly varying  with index
$s_0':=s/(\alpha_{1+}\vee\alpha_{1-})$. Note that
\begin{eqnarray*}
P(\sigma_1^2\geq x)&=& P(\eta_{0}>0)P\left(e^{u_1}\sigma_0^{2(\beta+\alpha_{1+})}\geq x\mid \eta_0>0\right)\\
&&+P(\eta_{0}<0)P\left(e^{u_1}\sigma_0^{2(\beta+\alpha_{1-})}\geq x\mid \eta_0<0\right).
\end{eqnarray*}
As $e^{u_1}$ admits regular variation of order $s'_0$, it admits a
moment of order $s'_0(\beta+\alpha_{1+})<s_0$. Note that
$E|\eta_0|^{2s}<\infty$ for any $s<s'$. An application of
Proposition \ref{prrvbis} thus gives
$E\left(\sigma_0^{2(\beta+\alpha_{1+})(s'_0+\iota)}\right)<\infty$
for $\iota>0$ small enough. By independence between $u_1$ and
$\sigma_0^2$ conditionally on $\eta_0>0$, we may apply a result by
Breiman (1965) to conclude that
$$P\left(e^{u_1}\sigma_0^{2(\beta+\alpha_{1+})}\geq x\mid
\eta_0>0\right)\sim
E\left(\sigma_0^{2(\beta+\alpha_{1+})s'_0}\right)P(e^u_1>x\mid
\eta_0>0),$$ as $x\to \infty.$ Applying the same arguments to
$P\left(e^{u_1}\sigma_0^{2(\beta+\alpha_{1-})}\geq x\mid
\eta_0<0\right) $ we obtain
\begin{eqnarray*}
P(\sigma_1^2\geq x)&\sim & P(\eta_{0}>0)E\left(\sigma_0^{2(\beta+\alpha_{1+})s'_0}\right)
P\left(e^{u_1}>x\mid \eta_0>0\right)\\
&&+P(\eta_{0}<0)E\left(\sigma_0^{2(\beta+\alpha_{1-})s'_0}\right)P\left(e^{u_1}>x\mid \eta_0<0\right)
\end{eqnarray*}
and the first assertion follows. The second assertion follows easily
by independence of $\eta_0$ and $\sigma_0$, with respective
regularly variation indexes $s'$ and $s'_0$. \zak

\subsection{Proof of Theorem \ref{consistency}} 
We will use the following standard result (see {\it  e.g.} Exercise
2.11 in Francq and Zakoian, 2010a).
\begin{lem}\label{petitlemme} Let $(X_n)$ be a sequence of random variables. If $\sup_n E|X_n|<\infty$, then almost surely
$n^{-1} X_n\to 0$ as $n\to\infty.$ The almost sure convergence may
fail when $\sup_n E|X_n|=\infty$. If the sequence $(X_n)$ is bounded
in probability, then $n^{-1} X_n\to 0$ in probability.
\end{lem}
Turning to the proof of Theorem~\ref{consistency}, first note that
{\bf A2}, {\bf A3} and Proposition~\ref{fracmoment} ensure the a.s.
absolute convergence of the series
\begin{equation}
\label{sigmadef}
\log \sigma_t^2(\btheta):={\cal B}_{\btheta}^{-1}(B)\left\{\omega+\sum_{i=1}^q\left(\alpha_{i+}1_{\{\epsilon_{t-i}>0\}}+\alpha_{i-}
1_{\{\epsilon_{t-i}<0\}}\right)\log\epsilon_{t-i}^2\right\}.
\end{equation}
Let \begin{equation} \label{critere} Q_n(\btheta)
=n^{-1}\sum_{t=r_0+1}^n\ell_t(\btheta),\qquad
 \ell_t(\btheta)= \frac{\epsilon_t^2}{\sigma_t^2(\btheta)} +\log \sigma_t^2(\btheta).
\end{equation}

Using standard arguments, as in the proof of Theorem 2.1 in Francq
and Zakoian (2004) (hereafter FZ),  the consistency
 is obtained by showing the following intermediate results
\begin{eqnarray*}
&&i) \; \lim_{n\to \infty} \sup_{\btheta\in\Theta}|
Q_n(\btheta)-\widetilde{Q}_n(\btheta)|=0\;\mbox{ a.s.; } \\
&& ii) \; \mbox{ if } \sigma_1^2(\btheta)=\sigma_1^2(\btheta_0) \mbox{ a.s. }\mbox{ then }
\btheta=\btheta_0;\\
&& iii) \; \mbox{ if } \btheta\ne \btheta_0\;, \quad
E\ell_t(\btheta)> E\ell_t(\btheta_0);\\
&& iv) \; \mbox{any $\btheta\neq \btheta_0$ has a neighborhood
$V(\btheta)$ such that }\\&&\qquad\qquad\qquad
\liminf_{n\to\infty}\inf_{\btheta^*\in
V(\btheta)}\widetilde{Q}_n(\btheta^*)>E\ell_t(\btheta_0)\;\; \mbox{ a.s. }
\end{eqnarray*}
 Because of the multiplicative form of the volatility, the step \emph{i)} is more delicate than in the standard GARCH case.
In the case $p=q=1$, we have
$$\log\sigma_t^2(\btheta)-\log\widetilde{\sigma}_t^2(\btheta)=\beta^{t-1}\left\{\log\sigma_1^2(\btheta)-\log\widetilde{\sigma}_1^2(\btheta)\right\},\quad\forall t\geq 1.$$
In the general case, as in FZ, using (\ref{condvi}) one can show
that for almost all trajectories,
\begin{equation}
\label{firstresult}
\sup_{\btheta\in\Theta}\left|\log\sigma_t^2(\btheta)-\log\widetilde{\sigma}_t^2(\btheta)\right|\leq K\rho^t,
\end{equation}
where $\rho\in(0,1)$ and $K>0$. 
First, we complete the
proof of $i)$  in the case  $p=q=1$ and $\alpha_+=\alpha_-$, for
which the notation is more explicit. In view of the multiplicative
form of the volatility
\begin{equation}
\label{multiform}
\sigma_t^2(\btheta)=e^{\beta^{t-1}\log\sigma_1^2(\btheta)}\prod_{i=0}^{t-2}e^{\beta^i\left\{\omega+\alpha\log\epsilon^2_{t-1-i}\right\}},
\end{equation}
we have
\begin{eqnarray*}\frac{1}{t}\log\left|\frac{1}{\sigma_t^2(\btheta)}-\frac{1}{\widetilde{\sigma}_t^2(\btheta)}\right|&=&\frac{-1}{t}\sum_{i=0}^{t-2}\beta^i\left\{\omega+\alpha\log\epsilon^2_{t-1-i}\right\}
\\&&+\frac{1}{t}\log\left|e^{-\beta^{t-1}\log\sigma_1^2(\btheta)}-e^{-\beta^{t-1}\log\widetilde{\sigma}_1^2(\btheta)}\right|.\end{eqnarray*}
Applying Lemma~\ref{petitlemme}, the first term of the right-hand
side of the equality tends almost surely to zero because it is
bounded by a variable of the form $|X_t|/t$, with $E|X_t|<\infty$,
under {\bf A5}. The second term is equal to
$$\frac{1}{t}\log\left|\left\{\log\sigma_1^2(\btheta)-\log\widetilde{\sigma}_1^2(\btheta)\right\}\beta^{t-1}e^{-\beta^{t-1}x^{*}}\right|,$$
where $x^*$ is between $\log\sigma_1^2(\btheta)$ and
$\log\widetilde{\sigma}_1^2(\btheta)$. This second term thus tends
to $\log|\beta|<0$ when $t\to\infty$. It follows that
\begin{equation}
\label{secondresult}
\sup_{\btheta\in\Theta}\left|\frac{1}{\sigma_t^2(\btheta)}-\frac{1}{\widetilde{\sigma}_t^2(\btheta)}\right|\leq K\rho^t,
\end{equation}
where $K$ and $\rho$ are as in (\ref{firstresult}). Now consider the
general case. Iterating (\ref{logGARCH}), using the compactness of
$\Theta$ and the second part of {\bf A2}, we have
\begin{eqnarray*}
\log \sigma_t^2(\btheta)&=&\sum_{i=1}^{t-1}c_i(\btheta)+c_{i+}(\btheta)1_{\{\epsilon_{t-i}>0\}}\log \epsilon_{t-i}^2+c_{i-}(\btheta)1_{\{\epsilon_{t-i}<0\}}\log \epsilon_{t-i}^2\\
&&+\sum_{j=1}^{p}c_{t,j}(\btheta)\log \sigma_{q+1-j}^2(\btheta)
\end{eqnarray*}
with
\begin{equation}
\label{infrhoi}
\sup_{\btheta\in\Theta}\max\{|c_i(\btheta)|,|c_{i+}(\btheta)|,|c_{i-}(\btheta)|,|c_{i,1}(\btheta)|,\dots,|c_{i,p}(\btheta)|\}\leq K\rho^i,\quad  \rho\in (0,1).
\end{equation}
We then obtain a multiplicative form for $\sigma_t^2(\btheta)$ which
generalizes (\ref{multiform}), and deduce that
$$\frac{1}{t}\log\left|\frac{1}{\sigma_t^2(\btheta)}-\frac{1}{\widetilde{\sigma}_t^2(\btheta)}\right|=a_1+a_2,$$
where
$$a_1=\frac{-1}{t}\sum_{i=1}^{t-1}c_i(\btheta)+c_{i+}(\btheta)1_{\{\epsilon_{t-i}>0\}}\log \epsilon_{t-i}^2+c_{i-}(\btheta)1_{\{\epsilon_{t-i}<0\}}\log \epsilon_{t-i}^2\to 0\quad\mbox{ a.s.}$$
in view of (\ref{infrhoi}) and Lemma~\ref{petitlemme}, and for
$x_j^*$'s between $\log \sigma_{q+1-j}^2(\btheta)$ and $\log
\widetilde{\sigma}_{q+1-j}^2(\btheta)$, {\small
\begin{eqnarray*}
a_2&=&\frac{1}{t}\log \left|\exp\left\{-\sum_{j=1}^{p}c_{t,j}(\btheta)\log \sigma_{q+1-j}^2(\btheta)\right\}-\exp\left\{-\sum_{j=1}^{p}c_{t,j}(\btheta)\log \widetilde{\sigma}_{q+1-j}^2(\btheta)\right\}\right|\\
&=&\frac{1}{t}\log \left|-\sum_{j=1}^{p}c_{t,j}(\btheta)\left\{\log \sigma_{q+1-j}^2(\btheta)-\log \sigma_{q+1-j}^2(\btheta)\right\}\exp\left\{-\sum_{k=1}^{p}c_{t,k}(\btheta)\log x_k^*\right\}\right|\\
&=&\frac{1}{t}\log \left|-\sum_{j=1}^{p}c_{t,j}(\btheta)\right|+o(1)\quad\mbox{  a.s.}
\end{eqnarray*}}
using  (\ref{condvi}) and (\ref{infrhoi}). Using again
(\ref{condvi}), it follows that $\limsup_{n\to\infty}a_2\leq
\log\rho<0$. We conclude that (\ref{secondresult}) holds true in the
general case. The proof of \emph{i)} then follows from
(\ref{firstresult})-(\ref{secondresult}), as in FZ.

To show \emph{ii)}, note that we have
\begin{equation}
\label{commentlappeler}
{\cal B}_{\btheta}(B)\log\sigma_t^2(\btheta)=\omega+
{\cal A}^+_{\btheta}(B)1_{\{\epsilon_{t}>0\}}\log\epsilon_{t}^2+
{\cal A}^-_{\btheta}(B)1_{\{\epsilon_{t}<0\}}\log\epsilon_{t}^2.
\end{equation}
 If $\log\sigma_1^2(\btheta)=\log\sigma_1^2(\btheta_0)$ a.s., by stationarity
we have $\log\sigma_t^2(\btheta)=\log\sigma_t^2(\btheta_0)$ for all
$t$, and thus we have almost surely
\begin{eqnarray*}&&\left\{\frac{{\cal A}^+_{\btheta}(B)}{{\cal
B}_{\btheta}(B)}-\frac{{\cal A}^+_{\btheta_0}(B)}{{\cal
B}_{\btheta_0}(B)}\right\}1_{\{\epsilon_{t}>0\}}\log\epsilon_{t}^2+
\left\{\frac{{\cal A}^-_{\btheta}(B)}{{\cal
B}_{\btheta}(B)}-\frac{{\cal A}^-_{\btheta_0}(B)}{{\cal
B}_{\btheta_0}(B)}\right\}1_{\{\epsilon_{t}<0\}}\log\epsilon_{t}^2\\&=&\frac{\omega_0}{{\cal
B}_{\btheta_0}(1)}-\frac{\omega}{{\cal B}_{\btheta}(1)}.\end{eqnarray*} Denote by
$R_t$ any random variable which is measurable with respect to
$\sigma\left(\{\eta_u,u\leq t\}\right)$. If
\begin{equation}
\label{hypopasvraie}
\frac{{\cal A}^+_{\btheta}(B)}{{\cal B}_{\btheta}(B)}\neq \frac{{\cal A}^+_{\btheta_0}(B)}{{\cal B}_{\btheta_0}(B)}\quad\mbox{ or }\quad
\frac{{\cal A}^-_{\btheta}(B)}{{\cal B}_{\btheta}(B)}\neq \frac{{\cal A}^-_{\btheta_0}(B)}{{\cal B}_{\btheta_0}(B)},
\end{equation}
 there exists a non null $(c_+,c_-)'\in \mathbb{R}^2$,
 such that $$c_+ 1_{\{\eta_t>0\}}\log\epsilon_t^2+c_- 1_{\{\eta_t<0\}}\log\epsilon_t^2+R_{t-1}=0\mbox{ a.s.}$$  This is equivalent to the two equations
$$\left(c_+ \log\eta_t^2+c_+ \log\sigma_t^2+R_{t-1}\right)1_{\{\eta_t>0\}}=0$$  and $$\left(c_- \log\eta_t^2+c_- \log\sigma_t^2+R_{t-1}\right)1_{\{\eta_t<0\}}=0.$$ Note that if an equation of the form $a\log x^2 1_{\{x>0\}}+b1_{\{x>0\}}=0$ admits two positive solutions then $a=0$. This result, {\bf A3}, and the independence between $\eta_t$ and $(\sigma_t^2,R_{t-1})$ imply that $c_+=0$. Similarly we obtain $c_-=0$, which leads to a contradiction. We conclude that (\ref{hypopasvraie}) cannot hold true, and the conclusion follows from {\bf A4}.

Since $\sigma_t^2(\btheta)$ is not bounded away from zero, the
beginning of the proof of \emph{iii)} slightly differs from that
given by FZ in the standard GARCH case. In view of
(\ref{commentlappeler}), the second part of {\bf A2} and {\bf A5}
entail that $E|\log\sigma_t^2(\btheta)|<\infty$ for all
$\btheta\in\Theta$. It follows that $E\ell_t^-(\btheta)<\infty$ and
$E|\ell_t(\btheta_0)|<\infty$.

The rest of the proof of \emph{iii)}, as well as that of \emph{iv)},
are identical to those given in FZ. \zak

\subsection{Proof of Theorem \ref{normality}} A Taylor
expansion gives
$$
\nabla_iQ_n(\widehat\btheta_n)-\nabla_i Q_n( \btheta_0)=\mathbb H_{i.} Q_n(\widetilde\btheta_{n,i})(\widehat\btheta_n-\btheta_0)\quad\mbox{for all}~1\le i\le d,
$$
where the $\widetilde\btheta_{n,i}$'s are such that
$\|\widetilde\btheta_{n,i}-\btheta_0\|\leq\|\widehat\btheta_n-\btheta_0\|$.
As in Section 5 of Bardet and Wintenberger (2009), the asymptotic
normality is obtained by showing:
\begin{enumerate}
\item  $n^{1/2}\nabla Q_n(\btheta_0)\to \mathcal N(\bzero,(\kappa_4-1)\bf J)$,
\item  $\|\mathbb H Q_n(\widetilde\btheta_n)-\bf J\|$ converges a.s.
to $0$ for any sequence $(\widetilde \btheta_n)$ converging a.s.
to $\btheta_0$ and  $ \bf J$  is invertible,
\item  $ n^{1/2}\|\nabla \widetilde Q_n(\widehat\btheta_n)-\nabla Q_n(\widehat\btheta_n)\|$
  converges a.s. to $0$.
\end{enumerate}
In order to prove the points 1-3 we will use the following Lemma
\begin{lem}\label{lem:mom}
Under  the assumptions of Theorem \ref{consistency} and {\bf A7},
for any $m>0$ there exists a neighborhood $\mathcal V $ of
$\btheta_0$ such that $E[\sup_{\mathcal V }
(\sigma_t^2/\sigma_t^2(\btheta))^m]<\infty$ and $E[\sup_{\mathcal V
} |\log\sigma_t^2(\btheta)|^m]<\infty$.
\end{lem}
\noindent{\bf Proof.} We have \begin{eqnarray*}
\log\sigma^2_t(\btheta_0)-\log\sigma^2_t(\btheta)&=&\omega_0-\omega+\sum_{j=1}^p\beta_j
\{\log\sigma^2_{t-j}(\btheta_0)-\log\sigma^2_{t-j}(\btheta)\}
\\&&\hspace{-3cm}+{\bf V}_{\btheta_0-\btheta}\bsigma_{t-1,r}+\mathcal A^+_{\btheta_0-\btheta}(B)1_{\eta_{t}>0}\log\eta^2_{t} +
\mathcal A^-_{\btheta_0-\btheta}(B)1_{\eta_{t}<0}\log\eta^2_{t}
\end{eqnarray*}
with $\bsigma_{t,r}=(\log \sigma_t^2(\btheta_0),\dots,
\log \sigma_{t-r+1}^2(\btheta_0))'$, 
$$
{\bf V}_{\btheta}=( \alpha_{1+} 1_{\{\eta_{t-1}>0\}}+ \alpha_{1-} 1_{\{\eta_{t-1}<0\}}+ \beta_{1},\ldots,
 \alpha_{r+} 1_{\{\eta_{t-r}>0\}}+ \alpha_{r-} 1_{\{\eta_{t-r}<0\}}+ \beta_{r}).
$$
Under {\bf A2}, we then have
\begin{multline*}
\log\sigma^2_t(\btheta_0)-\log\sigma^2_t(\btheta)=\mathcal B_{\btheta}^{-1}(B)\left\{\omega_0-\omega+{\bf V}_{\btheta_0-\btheta}\bsigma_{t-1,r}\right.\\
\left.+(\mathcal
A^+_{\btheta_0-\btheta}(B)1_{\eta_{t}>0}\log\eta^2_{t}+\mathcal
A^-_{\btheta_0-\btheta}(B)1_{\eta_{t}<0}\log\eta^2_{t}\right\}.
\end{multline*}
Under {\bf A7} the assumptions of Proposition~\ref{prmom} hold. From
the proof of that proposition, we thus have that $E\exp(\delta
\|\mathrm{Abs}(\bsigma_{t,r})\|)$ is finite for some $\delta>0$.

Now, note that ${\bf V}_{\btheta}$, $\mathcal A_{\btheta}^+(1)$ and
$\mathcal A_{\btheta}^+(1)$ are continuous functions of $\btheta$.
Choosing a sufficiently small neighborhood $\mathcal V$ of
$\btheta_0$,  one can  make $\sup_{\mathcal V}\|{\bf
V}_{\btheta_0-\btheta}\|$, $\sup_{\mathcal V}|\mathcal
A_{\btheta_0-\btheta}^+(1)|$ and $\sup_{\mathcal V}|\mathcal
A_{\btheta_0-\btheta}^+(1)|$ arbitrarily small. Thus
$E[\exp(m\sup_{\mathcal V}\|{\bf
V}_{\btheta_0-\btheta}{\bsigma}_{t,r}\|)]$ and
$E[\exp(m\sup_{\mathcal V}\|(\mathcal
A^+_{\btheta_0-\btheta}(B)1_{\eta_{t-1}>0} +\mathcal
A^-_{\btheta_0-\btheta}(B)1_{\eta_{t-1}<0})\log(\eta^2_{t-1})\|)]$
are finite for an appropriate choice of $\mathcal V$ depending on $m$. 
We conclude that $E\left[\exp\left(m\sup_{\mathcal
V}\left|\log\left\{\sigma^2_t(\btheta_0)/\sigma^2_t(\btheta)\right\}\right|\right)\right]<\infty$
and the first assertion of the lemma is proved.

Consider now the second assertion. We have
$$\sup_{\mathcal V}| \log\sigma_t^2(\btheta)|\le |\log \sigma_t^2|+ \sup_{\mathcal V}|\log(\sigma_t^2(\btheta_0)/\sigma_t^2(\btheta))|. $$
We have already shown that the second term admits a finite moment of
order $m$. So does the first term, under {\bf A7}, by Proposition
\ref{anylogmoment}.$\hfill\square$

Now let us prove  the point 1.  
In view of (\ref{critere}) we have
$$
\nabla Q_n(\btheta)=\frac1n\sum_{t=r_0+1}^n \left(1-\frac{\epsilon^2_t}{\sigma_t^2(\theta)}\right)
\nabla \log\sigma_t^2(\btheta) \quad \mbox{and thus}
\quad \nabla Q_n(\btheta_0)=\frac1n\sum_{t=r_0+1}^n (1-\eta^2_t)\nabla \log\sigma_t^2(\btheta_0).
$$
Because $\eta_t$ and $\log\sigma_t^2(\btheta_0)$ are independent,
and since $E \eta_t^2=1$, the Central Limit Theorem for martingale
differences applies (see Billingsley (1961)) whenever
$\bQ=(\kappa_4-1)E(\nabla \log\sigma_t^2(\btheta_0)\nabla
\log\sigma_t^2(\btheta_0)')$ exists. For any
$\btheta\in\stackrel{\circ}{\Theta}$, the random vector $\nabla
\log\sigma_t^2(\btheta)$ is the stationary solution of the equation
\begin{equation}\label{eq:nabla}
\nabla \log\sigma_t^2(\btheta)=\sum_{j=1}^p\beta_j \nabla \log\sigma^2_{t-j}(\btheta)+ \left(\begin{array}{c}1 \\\bepsilon_{t-1,q}^+\\
\bepsilon_{t-1,q}^-\\\bsigma^2_{t-1,p}(\btheta)\end{array}\right),
\end{equation}
where $ \bsigma^2_{t,p}(\btheta)=(\log \sigma_{t}^2(\btheta),
\ldots, \log \sigma_{t-p+1}^2(\btheta))'$.

Assumption {\bf A2} entails that $\nabla \log\sigma_t^2(\btheta)$ is
a linear combination of $\bepsilon_{t-i,q}^+$, $\bepsilon_{t-i,q}^-$
and $\log \sigma_{t-i}^2(\btheta)$ for $i\geq 1$. Lemma
\ref{lem:mom} ensures that, for any $m>0$, there exists a
neighborhood $\mathcal V$ of $\btheta_0$ such that $E[\sup_{\mathcal
V}|\log \sigma_{t-i}^2(\btheta)|^m]<\infty$. By
Proposition~\ref{anylogmoment},  $\bepsilon_{t-i,q}^+$ and
$\bepsilon_{t-i,q}^-$ admit moments of any order. Thus, for any
$m>0$ there exists $\mathcal V$ such that $E[\sup_{\mathcal
V}\|\nabla \log\sigma_t^2(\btheta)\|^m]<\infty$. In particular,
$\nabla \log\sigma_t^2(\btheta_0)$  admits moments of any order.
Thus point 1. is proved.

Turning to point 2., we have
$$
\mathbb H Q_n(\btheta) =n^{-1}\sum_{t=r_0+1}^n\mathbb H \ell_t(\btheta),
$$
where
\begin{equation}\label{eq:l}
\mathbb H \ell_t(\btheta)=\left(1-\frac{\eta_t^2\sigma_t^2(\btheta_0)}{\sigma_t^2(\btheta)}\right)\mathbb H \log
\sigma_t^2(\btheta)+\frac{\eta_t^2\sigma_t^2(\btheta_0)}{\sigma_t^2(\btheta)}\nabla
\log \sigma_t^2(\btheta) \nabla \log
\sigma_t^2(\btheta)'.\end{equation}
By Lemma \ref{lem:mom}, the term
$\sigma_t^2(\btheta_0)/\sigma_t^2(\btheta)$ admits moments of order
as large as we need uniformly on a well chosen neighborhood
$\mathcal V$ of $\btheta_0$. Let us prove that it is also the case
for $\mathbb H \log\sigma_t^2(\btheta)$. Computation gives
\begin{eqnarray*}
\mathbb H \log\sigma_t^2(\btheta)=\sum_{j=1}^p\beta_j \mathbb H \log\sigma_{t-j}^2(\btheta)+
\left(\begin{array}{c}\bzero_{(2q+1)\times d}\\\nabla'\bsigma^2_{t-1,p}(\btheta) \end{array}\right)+\left(\begin{array}{c}\bzero_{(2q+1)\times d}\\\nabla'\bsigma^2_{t-1,p}(\btheta) \end{array}\right)'.
\end{eqnarray*}
From this relation and  {\bf A2} we obtain
$$
\mathbb H \log\sigma_t^2(\btheta)=
\left(\begin{array}{c}\bzero_{(2q+1)\times d}\\\mathcal B_{\btheta}(B)^{-1}\nabla'\bsigma^2_{t-1,p}(\btheta) \end{array}\right)+
\left(\begin{array}{c}\bzero_{(2q+1)\times d}\\\mathcal B_{\btheta}(B)^{-1}\nabla'\bsigma^2_{t-1,p}(\btheta) \end{array}\right)'.
$$
Thus $\mathbb H \log\sigma_t^2(\btheta)$ belongs to $\mathcal
C(\mathcal V)$ and  is integrable  because we can always choose
$\mathcal V$ such that $\sup_{\mathcal
V}\|\nabla'\bsigma^2_{t-1,p}(\btheta)\|\in L^m$ (see the proof of
point 1. above).

An application of the Cauchy-Schwarz inequality  in the RHS term of
\eqref{eq:l} yields  the integrability of $\sup_{\mathcal V}\mathbb
H \ell_t(\btheta)$. The first assertion of point 2. is proved by an
application of the ergodic theorem on $(\mathbb H \ell_t(\btheta))$
in the Banach space $\mathcal C(\mathcal V)$ equipped with the
supremum norm:
$$
\sup_{\mathcal V}\|\mathbb H Q_n(\btheta)-E[\mathbb H \ell_0(\btheta)]\| \to 0\qquad a.s.
$$
An application of Theorem \ref{consistency} ensures that
$\hat\btheta_n$ belongs a.s. to $\mathcal V$ for sufficiently large
$n$. Thus
$$
\|\mathbb H Q_n(\hat\btheta_n)-E[\mathbb H \ell_0(\btheta_0)]\|\le \sup_{\mathcal V}\|\mathbb H Q_n(\btheta)-E[\mathbb H \ell_0(\btheta)]\| +\|E[\mathbb H \ell_0(\hat\btheta_n)]-E[\mathbb H \ell_0(\btheta_0)]\|
$$
converges a.s. to $0$ by continuity of $\btheta\to E[\mathbb H
\ell_0( \btheta)]$ at $\btheta_0$ as a consequence of a dominating
argument on $\mathcal V$. The first assertion of point 2. is proved.
The invertibility of matrix $\bJ$ follows from arguments used in the
proof of Theorem \ref{consistency}, ii).
\\


From \eqref{eq:nabla} and an equivalent representation for $\nabla
\log\widetilde \sigma_t^2( \btheta )$, we have
\begin{eqnarray*}
\nabla \log
\sigma_t^2( \btheta )-\nabla \log\widetilde
\sigma_t^2( \btheta )&=&\sum_{j=1}^p\beta_j (\nabla \log\sigma_{t-j}^2( \btheta )-\nabla \log\widetilde
\sigma_{t-j}^2( \btheta ))\\
&&+ \left(\begin{array}{c} \bzero_{2q+1}\\ \bsigma^2_{t-1,p}(\btheta)-\widetilde \bsigma^2_{t-1,p}(\btheta) \end{array}\right)
\end{eqnarray*}
where $\widetilde \bsigma^2_{t,p}$ is defined as $\bsigma^2_{t,p}$.
Thus, there exist continuous functions $d_i$ and $d_{t,i}$ defined
on $\Theta$ such that
\begin{eqnarray*}
\nabla \log \sigma_t^2( \btheta )-\nabla \log\widetilde
\sigma_t^2( \btheta )&=&\sum_{i=1}^{t-1}d_i(\btheta)(\log\sigma^2_{t-i}(\btheta)-\log\widetilde\sigma^2_{t-i}(\btheta))\\
&&+\sum_{j=1}^pd_{t,j}(\btheta)\nabla \log\sigma_{p+1-j}^2(\btheta).
\end{eqnarray*}
The sequences of functions $(d_i)$, $(d_{i,j})$, $1\le j\le p$,
satisfy the same uniform rate of convergence as the functions $c_i$,
$c_{i+}$, $c_{1-}$ and $c_{i,j}$ in \eqref{infrhoi}. An application
of \eqref{firstresult} yields the existence of $K>0$ and
$\rho\in(0,1)$ such that $\sup_\Theta \|\nabla \log \sigma_t^2(
\btheta )-\nabla \log\widetilde \sigma_t^2( \btheta )\|\le K\rho^t$,
for almost all trajectories. Point 3.  easily follows and the
asymptotic normality is proved.\zak

\section{Conclusion}
\label{sec8}

In this paper, we investigated the probabilistic properties of the
log-GARCH($p,q$) model. We found sufficient conditions for the
existence of moments and log-moments of the strictly stationary
solutions. We analyzed the dependence structure through the
 leverage effect and the regular
variation properties, and we compared this structure with that of
the EGARCH model.

As far as the estimation is concerned, it should be emphasized that
the log-GARCH model appears to be much more tractable than the
EGARCH. Indeed, we established the strong consistency and the
asymptotic normality of the QMLE under mild assumptions. For EGARCH
models, such properties have only been established for the
first-order model and with strong invertibility constraints (see
Wintenberger, 2013). By comparison with standard GARCH, the
log-GARCH model is not more difficult to handle: on the one hand,
the fact that the volatility is not bounded below requires an
additional log-moment assumption, but on the other hand the
parameters are nor positively constrained.


A natural extension of this work,
 aiming at pursuing the
comparison between the two classes of models, would rely on
 statistical tests. By embedding the
log-GARCH model in a more general framework including the log-GARCH,
it should be possible to consider a LM test of the log-GARCH null
assumption. 
Another problem of interest would be to check  validity of the
estimated models.
We leave these issues for further investigation, viewing the results
of this paper as a first step in these directions.

\bigskip
\begin{center}
{\Large\bf References}
\end{center}

\begin{description}
\item[Allen, D. ,  Chan, F.,  McAleer, M.,  and  Peiris, S.] (2008) Finite sample properties of the QMLE for the log-ACD model: application to Australian stocks. {\it  Journal of Econometrics} 147, 163--183.
\item[Bardet, J.-M. and Wintenberger, O.] (2009)
Asymptotic normality of the Quasi Maximum Likelihood estimator
for multidimensional causal processes. {\it Annals of Statistics} 37,
2730--2759.
\item[Bauwens, L. and  Giot, P.] (2000) The logarithmic ACD model: An application to the bidask
quote process of three NYSE stocks. {\it  Annales D'Economie et
de Statistique} 60, 117--145.
\item[Bauwens, L. and  Giot, P.] (2003) Asymmetric ACD models: introducing price information in ACD models.
{\it  Empirical Economics} 28, 709--731.
\item[Bauwens, L., Galli, F. and Giot, P.] (2008) The moments of log-ACD models. {\it  QASS} 2, 1--28.
\item[Bauwens, L., Giot, P., Grammig, J. and  Veredas, D.] (2004) A comparison of financial duration models via density forecast. {\it  International Journal of Forecasting} 20, 589--604.
\item[Berkes, I., Horv\'ath, L. and  Kokoszka, P.] (2003)
GARCH processes: structure and estimation. {\it Bernoulli} 9,
201--227. \item[Billingsley, P.](1961) The Lindeberg-Levy
theorem for martingales. {\it  Proceedings of the American
Mathematical Society} 12,  788--792.
\item[Bollerslev, T.] (1986)
Generalized autoregressive conditional heteroskedasticity. {\it
Journal of Econometrics} 31, 307--327.
\item[Bougerol, P. and   Picard, N.] (1992a)
Strict stationarity of generalized autoregressive processes. {\it
Annals of Probability} 20, 1714--1729.
\item[Bougerol, P. and Picard, N.] (1992b)
Stationarity of GARCH processes and of some nonnegative time
series. {\it  Journal of Econometrics} 52, 115--127.
\item[Breiman, L.] (1965)
On some limit theorems similar to the arc-sin law. {\it Theory
of Probability and Applications} 10, 323--331.
\item[Drost, F.C.  and  Nijman, T.E.] (1993) Temporal aggregation of GARCH
 processes. {\it Econometrica} 61, 909--927.
\item[Engle, R.F.] (1982) Autoregressive conditional heteroskedasticity with estimates of the variance
of U.K. inflation.  {\it Econometrica} 50, 987--1008.
\item[Engle, R.F. and  Russell, J.R.] (1998) Autoregressive conditional duration: a new model for
irregularly spaced transaction data. {\it  Econometrica} 66,
1127--1162.
\item[Francq, C. and Zakoïan, J-M.] (2004) Maximum likelihood estimation of pure GARCH and ARMA-GARCH.  {\it Bernoulli} 10, 605--637.
\item[Francq, C. and Zakoïan, J-M.] (2010a) {\it GARCH Models : structure, statistical inference
and financial applications}. John Wiley.
\item[Francq, C. and
Zakoïan, J-M.] (2010b) Inconsistency of the MLE and inference
based on weighted LS for LARCH models. {\it Journal of
Econometrics} 159, 151--165.
\item[Geweke, J.] (1986) Modeling the persistence of conditional variances: a comment. {\it Econometric
Review} 5, 57--61.
\item[He, C., Ter\"{a}svirta, T. and Malmsten, H.] (2002) Moment structure of a family of first-order exponential GARCH models. {\it Econometric Theory} 18, 868--885.
\item[Kristensen, D. and A. Rahbek] (2009)
Asymptotics of the QMLE for non-linear ARCH models. {\it Journal
of Time Series Econometrics} 1,  Paper 2.
\item[Ling, S. and M. McAleer] (2002) Necessary
and sufficient moment conditions for the GARCH($r,s$) and
 asymmetric GARCH($r,s$) models. {\it Econometric Theory} 18,
 722--729.
\item[Mikosch,
T. and Rezapur, M.] (2012) Stochastic volatility models with possible extremal clustering. {\it Bernoulli} To appear.
\item[Mikosch, T. and Starica, C.] (2000) Limit theory for the sample autocorrelations and extremes of a GARCH (1,1) process. {\it Annals of Statistics} 28, 1427--1451.
\item[Milhøj, A.] (1987) A multiplicative parameterization of ARCH Models.  {\it Working paper}, Department
of Statistics, University of Copenhagen.
\item[Nelson D.B.] (1991)
Conditional heteroskedasticity in asset returns : a new
approach. {\it Econometrica} 59, 347--370.
\item[Pantula, S.G.] (1986) Modeling the persistence of conditional variances: a comment. {\it Econometric
Review} 5, 71--74.
\item[Robinson, P.M.] (1991). Testing for strong serial correlation and dynamic conditional heteroskedasticity
in multiple regression. {\it  Journal of Econometrics} 47,
67--84.
\item[Sucarrat, G. and Escribano, A.] (2010) The Power Log-GARCH Model. {\it Working document}, Economic Series 10-13, University Carlos III, Madrid.
\item[Sucarrat, G. and Escribano, A.] (2012) Automated model selection in finance: general-to-specific modelling of the mean and volatility specifications.
{\it Oxford Bulletin Of Economics And Statistics} 74, 716--735.
\item[Tweedie, R.L.] (2008)
Invariant measures for Markov chains with no irreducibility
Assumptions. {\it Journal of Applied Probability} 25, 275--285.
\item[Wintenberger, O.] (2013)
 Continuous invertibility and stable QML estimation of the
EGARCH(1,1) model. Preprint arXiv:1211.3292.
\end{description}

\end{document}